\def \Ext{\mathop{\rm Ext}\nolimits}
\def \la {\lambda}

\def \tl {{\tilde {\lambda}}}
\def \tm {{\tilde {\mu}}}
\def \arr {\rightarrow}
\def \larr {\longrightarrow}
\def \Wedge {\Lambda}

\magnification=1200

\centerline {\bf On the Ext groups between Weyl modules for ${\rm GL}_n$}

\vskip 1cm

\centerline {\bf Upendra Kulkarni}

\vskip .3cm

\centerline {Truman State University, Kirksville, MO 63501 (Email: kulkarni@truman.edu)}

\vskip 1cm

\centerline{\bf Abstract}

\vskip .5cm

This paper studies extension groups between certain Weyl modules for
the algebraic group ${\rm GL}_n$ over the integers. Main results include:
(1) A complete determination of Ext groups between Weyl modules whose
highest weights differ by a single root and (2) Determination of Ext$^1$
between an exterior power of the defining representation and any Weyl
module. The significance of these results for modular representation
theory of ${\rm GL}_n$ is discussed in several Remarks. Notably the first
result leads to a calculation of Ext groups between neighboring Weyl
modules for ${\rm GL}_n$ and also recovers the ${\rm GL}_n$ case of a
recent result of Andersen. Some generalities about Ext groups between
Weyl modules and a brief overview of known results about these groups
are also included.


\vskip 1cm

\centerline {INTRODUCTION}

\vskip .5cm

This paper studies some homological aspects of the representation
theory of the algebraic group ${\rm GL}_n$ over the integers. More
specifically we will investigate Ext groups between certain pairs of
Weyl modules for ${\rm GL}_n({\bf Z})$ and also discuss the significance of
these results in modular representation theory of ${\rm GL}_n$.
Weyl modules are universal highest weight modules in the representation
category of split reductive algebraic groups. Extensions between these modules
are of interest in representation theory. Let us outline and briefly
discuss the main results. For a highest weight $\la$, let $K_{\la}$
denote the corresponding Weyl module.

\vskip .3cm

(1) Theorem 2.1 gives a complete determination of the groups
Ext$^{i}_{{\rm GL}_n({\bf Z})}(K_{\la}, K_{\mu})$ where ${\mu}-{\la}$ is a
positive root $\alpha$ of ${\rm GL}_n$. Here Ext$^1$ is cyclic of order
${\langle} {\la} + {\rho} , {\alpha} \, \check{} \, {\rangle} + 1$,
where $\rho$ is half the sum of positive roots and all other Ext groups vanish.
By the Universal Coefficient Theorem, one then gets all the modular Ext groups
as well for this class of examples. Moreover, in the modular case and also over
the $p$-adic integers, one can calculate Ext groups between any two neighboring
Weyl modules for ${\rm GL}_n$ using Theorem 2.1 along with the translation principle.
This recovers the ${\rm GL}_n$ case of a recent result of [Andersen3]
obtained independently around the same time as an equivalent version of Theorem
2.1. Neighboring Weyl modules are defined only if regular weights exist, which
happens only if the characteristic $p$ is at least $n$. Theorem 2.1 is of additional
interest because it gives results in the modular case even for small primes.
See the Remarks after the proof of Theorem 2.1 for a discussion. In view of
Andersen's result, it seems natural to hope that Theorem 2.1, as stated at
the beginning of this paragraph, remains true for all split reductive algebraic
groups over {\bf Z}.

\vskip .3cm

(2) Theorem 2.2 determines Ext$^1$ between an exterior power of the defining
representation and any Weyl module $K_\la$. This group is cyclic and its order
is the {\it gcd} of several integers that can be described explicitly in terms
of the weight $\la$. By contravariant duality and conjugate symmetry of Ext groups,
this result also leads to the calculation of certain other Ext$^1$ groups, e.g., Ext$^1$
between a symmetric power of the defining representation and a dual Weyl module.
It is interesting to compare Theorem 2.2 with the known enumeration of the
composition factors of symmetric powers in characteristic $p$. See Remark 3
after the proof of Theorem 2.2.

\vskip .3cm

Let us review the previously known information about the groups
Ext$^i(K_\la, K_\mu)$, in the modular as well as the integral
setting. (Note that the Ext is always taken in the appropriate category
of representations. But, as in the previous sentence, the notation will
often omit this fact relying on the context to convey the intended
meaning.) Even though the work in this paper deals directly only
with the integral case, the two cases are intimately connected
via the Universal Coefficient Theorem.

\vskip .3cm

Let us begin with the modular case, i.e., that of reductive algebraic groups
over an algebraically closed field of characteristic $p >0$. An important result
here is the vanishing theorem of Cline-Parshall-Scott-van der Kallen. [CPSvdK]
proves that for any such group Ext$^i(K_\la, K_\mu)$ vanishes unless
$\la \leq \mu$ under the dominance partial order, i.e., unless $\mu - \la$
is a sum of positive roots. This vanishing result was strengthened after
the proof of the strong linkage principle (see, e.g., [Jantzen]). As a result,
one can replace the dominance relation $\leq$ in the vanishing theorem
above by a finer relation $\uparrow$, where $\uparrow$ is the ``linkage"
relation defined using the dot action of the affine Weyl group $W_p$ on
weights. When $\la \uparrow \mu$, the calculation of the groups
Ext$^i(K_\la, K_\mu)$ is not known in general and is likely to be difficult.
But we do know the answer in the important case of  ``neighboring Weyl
modules."  This refers to the siutation when $\la$ and $\mu$ are regular
(i.e., have a trivial stabilizer in $W_p$) and $\la  <  s.\la = \mu$, where
$s$ is the reflection in a wall of the alcove containing $\la$. In this case
[Jantzen, II.7.19] proves that Ext$^i(K_\la, K_{s.\la})$ is one-dimesional
if $i$ = 0 or 1 and vanishes otherwise. (The relevance of this result for us
will be discussed in the Remarks after the proof of Theorem 2.1.) Beyond this
there are several results in more or less special cases, in which the answers
as well as the needed arguments are often involved. See [Jantzen, II.6.25]
for a discussion of results regarding homomorphisms between Weyl modules and
[Wen] for some further cases. [O-M] calculates Hom between certain hook
representations for ${\rm GL}_n$. [Erdmann] and [CE] respectively calculate
Ext$^1$ and Ext$^2$ between modular Weyl modules for ${\rm SL}_2$.

\vskip .3cm

Before turning to the integral case, let us digress to comment on the case of
the Bernstein-Gelfand-Gelfand category $\cal \char'117 $ of representations
of complex semisimple Lie algebras. This situation is somewhat parallel to
modular representations of reductive algebraic groups. The universal highest
weight modules here (i.e., analogues of Weyl modules) are Verma modules.
The extent of our knowledge about the Ext groups between these is similar to
that for Weyl modules. Similar vanishing properties hold, and Hom groups
between neighboring Verma modules are one-dimensional. But additionally,
unlike in the case of Weyl modules, there is a well-known calculation of Hom
groups between arbitrary Verma modules. The answer involves the Bruhat
order on the associated Weyl group. A guess was made in [GJ] expressing all
Ext groups between arbitrary Verma modules in terms of $R$-polynomials for
the associated Weyl group. But this guess was found to be incorrect in [Boe],
underscoring the seeming difficulty of calculating these Ext groups.

\vskip .3cm

Let us now survey the integral case, i.e., that of Ext groups between integral
Weyl modules. First of all, the vanishing result of [CPSvdK] immediately
carries over to this case.
The strengthening due to the linkage principle is not available over {\bf Z}.
But the Universal Coefficient Theorem does allow us to translate the
characteristic $p$ results mentioned above into information about $p$-torsion
of integral Ext groups between Weyl modules. So we know that there will be
no $p$-torsion in Ext$^{i}_{{\rm GL}_n({\bf Z})}(K_{\la}, K_{\mu})$ unless
$\la \uparrow \mu$ under the action of $W_p$. Similarly a nonzero
homomorphism (respectively, a one-dimensional Hom group) between two
Weyl modules in characteristic $p$ translates into nonvanishing (respectively
cyclic) $p$-torsion in the corresponding integral Ext$^1$. Beyond this what
we have are mainly special case results for the Ext groups between Weyl
modules for ${\rm GL}_n$ due to various authors. These results are obtained
following the basic approach in [AB2], which relies on constructing explicit
projective resolutions of Weyl modules. See [AB2], [F], [BF] for Ext$^1$
between special pairs of representations of ${\rm GL}_2$ and ${\rm GL}_3$.
See [Akin] and [Maliakas] for Ext$^1$ in some other special cases involving
hooks. (The calculation in [Maliakas] implies the earlier modular result in [O-M].)
See [R-G] for Ext$^2$ for ${\rm GL}_2$. (The results of [AB2] and [R-G] together
recover the modular calculation of Ext$^1$ for ${\rm SL}_2$ cited above. [Erdmann,
p. 456] describes this calculation in terms of a set $\Psi(r)$ which the author finds
convenient to define as a union of two subsets. It turns out that via the Universal
Coefficient Theorem, the first of these subsets is accounted for by integral Ext$^2$
[R-G] and the second by integral Ext$^1$ [AB2].) Theorem 2.1 in this paper
generalizes a result of [Maliakas] and intersects to various degrees with the results
of other papers. Theorem 2.2 generalizes certain Ext$^1$ calculations in [AB2] and
[Akin].

\vskip .3cm

Unlike [AB2], the approach here does not use resolutions directly. Instead a
key tool will be the Skew Representative Theorem from [Kulkarni1]. Some
other notions and methods that we will use are as follows. (These topics are
discussed in Section 1.) The explicit combinatorial/multilinear-algebraic
descriptions of Weyl modules for ${\rm GL}_n$, ordinary as well as ``skew," due
to Akin-Buchsbaum-Weyman will be very useful for us in proofs. Also very
useful will be Pieri-type rules giving filtrations of certain special skew Weyl
modules. The Schur algebras $S(n,r)$ will play a role, though mostly in the
background. These algebras were first treated systematically in [Green],
where it is proved that polynomial representations of ${\rm GL}_n$ of degree $r$
are equivalent to the representations of the Schur algebra $S(n,r)$. Then the
main new idea can be described as a way to reduce questions about a Schur
algebra to questions about another Schur algebra of smaller degree and then
to use recursion. The Skew Representative Theorem is the vehicle that allows
one to do this.

\vskip .3cm

Finally, let me indicate why this paper deals only with ${\rm GL}_n$ and
comment on possible generalizations. (Incidentally the main results here
stay valid for ${\rm SL}_n$. We will use ${\rm GL}_n$ as it will be more natural to
for us work in that setting.)
Many of the generalities used in this paper for ${\rm GL}_n$ hold for other
reductive algebraic groups too. Donkin has defined algebras generalizing
$S(n,r)$ for all split reductive algebraic groups, now known as Schur algebras.
These more general Schur algebras are examples of the quasihereditary
algebras of Cline-Parshall-Scott. Several of the properties relevant to this
paper (existence of suitable filtrations, triangular Ext-vanishing properties
with respect to a suitable partial order on weights) hold in the broader
setting of highest weight categories (i.e., representations of quasihereditary
algebras) of Cline-Parshall-Scott or BGG categories of Irving. But after the
generalities, what allows one to push through with the calculations is the
availability for ${\rm GL}_n$ of some very explicit characteristic-free constructions.
To carry out a similar method for other reductive groups, one should look
for a good analogue of skew Weyl modules for which a version of the Skew
Representative Theorem holds. Donkin has previously constructed ``skew
modules" for reductive groups [Donkin2].
More recently, he has also proved [Donkin3] a version of the Skew
Representative Theorem for these general skew modules. Since all his
constructions are abstract, it is not immediately clear how one can
use them to get analogues of results for ${\rm GL}_n$ that are obtained in
this paper by concrete calculations.

\vskip .3cm

In another direction, natural analogues of all the results in this paper
should hold for quantum ${\rm GL}_n$ as well. One just has to replace the
[ABW] constructions by the constructions of Hashimoto-Hayashi for
quantum Schur and Weyl modules. But we will not pursue this here.

\vskip .3cm

Since the first version of this paper was written, a result similar to
Theorem 2.1 and obtained independently around the same time has been published
by Andersen. For any reductive algebraic group over an algebraically closed
field of characteristic $p>0$, [Andersen3] calculates Ext groups between
neighboring Weyl modules for the corresponding Chevalley group over the
$p$-adic integers.  The connection between Theorem 2.1 and Andersen's result
is discussed after the proof of Theorem 2.1. One upshot is that it seems
natural to hope for the validity of Theorem 2.1 for any split reductive
algebraic group over {\bf Z}. A possible way to approach the expected
generalization is to use Donkin's work mentioned above. Another possibility
is to try to extend Andersen's proof, namely work with one prime at a time
and use translation functors, keeping careful track of the modules that arise
for small primes.

\vskip 1cm

\centerline {1. BACKGROUND AND NOTATION}

\vskip .5cm

This section discusses, in a little more detail than is strictly
necessary, the following two topics. First, the definition and
some properties of Weyl and Schur modules (both ordinary and
``skew") following Akin-Buchsbaum-Weyman, and second, some results
of a general nature about certain Ext groups of interest.
Proofs are given for a couple of easy results apparently
unrecorded elsewhere. Along the way we will also discuss the role
of the Schur algebras $S(n,r)$ and then use these algebras while
discussing Ext groups.

\vskip .3cm

Throughout the rest of this paper, the ground ring will be {\bf Z}
unless otherwise indicated. The exceptions will mainly occur in some
of the Remarks after the results, where the significance of the results
in modular representation theory is discussed. In any event all
discussions where we will need to consider a ground ring other than
{\bf Z} will take place in a separate paragraph, with an appropriate
notice to that effect at the beginning of that paragraph.

\vskip .3cm

We will study some homological aspects of the representation theory of
the reductive algebraic group scheme ${\rm GL}(F)={\rm GL}_n$, where $F$ is a free
abelian group of rank $n$. See, e.g., [Jantzen]. Weights for ${\rm GL}_n$ are
multiplicative characters of a maximal toral subgroup scheme. We will
follow the common practice of taking this to be the diagonal subgroup
$Diag({\bf Z}^n)$ and identifying weights for ${\rm GL}_n$ with $n$-tuples of
integers $\la = (\la_1, \ldots, \la_n)$. A weight $\la$ is a polynomial
weight if $\la_i$ are all nonnegative and then the degree of $\la$ is
$|\la| = \sum {\la_i}$. Dominant weights are the ones with
$\la_1 \geq \ldots \geq \la_n$. Dominant polynomial weights are thus just
partitions (with at most $n$ parts), which we will frequently identify
with their Young diagrams (with at most $n$ rows). The usual dominance
partial order on weights is the one generated by stipulating that
$( \ldots, \la_i, \la_{i+1}, \ldots) < (\ldots ,\la_i+1 , \la_{i+1}-1, \ldots)$.
For Young diagrams this means that moving boxes upwards gives a bigger
partition. Note that a given Young diagram $\la$ can be considered a
weight for all ${\rm GL}_n$ with $n \geq$ the number of rows of $\la$. In the
main results it will be convenient for us to adopt this point of view by
fixing the partitions and dealing with all such $n$ simultaneously.

\vskip .3cm

Given any partition $\la$, its conjugate $\tl$ is the partition whose
diagram is obtained by transposing (i.e., by exchanging rows and columns
of) the diagram of $\la$. In the literature the notation $\la^{\prime}$
is often used instead of $\tl$, but we will follow the notation in [ABW].
It should be noted that even if $\la$ has at most $n$ nonzero parts (and
so is a weight for ${\rm GL}_n$), $\tl$ may have more nonzero parts and so need
not be a weight for ${\rm GL}_n$. Thus conjugation is really an operation on
``stable" weights.

\vskip .3cm

{\bf Weyl and Schur modules.}
Let $\la$ be a partition having number of rows $\leq$ the rank of $F$.
Following [ABW], we will use $K_\la(F)$ to denote the Weyl module of
highest weight $\la$. When there is no danger of confusion, we will often
drop the $F$ and simply write $K_\la$. {\it The notation in this paper
is different from the standard one!} The standard notation for the Weyl
module of highest weight $\la$ is $V(\la)$. See, e.g., [Jantzen]. But we
will need to use heavily the constructions in [ABW] of Weyl modules as
well as of certain generalizations of Weyl modules defined in the same
paper (see below). So we will sacrifice the standard notation for the
sake of consistency with this reference. Note further that in [ABW] Weyl
modules are called coSchur modules, a term that will not be used in this
paper (and which was abandoned later by its authors as well).

\vskip .3cm

Some remarks are in order before introducing the ``skew Weyl modules."
First of all, as defined in [ABW], $K_\la$ is a functor that assigns to
any finitely generated free abelian group $F$ an abelian group $K_\la(F)$.
That $K_\la$ is a functor means in particular that $K_\la(F)$ is a
representation of ${\rm GL}(F)$. Also note that the [ABW] definition of
$K_\la$ works even when the number of rows in $\la$ is greater than
the rank of $F$, i.e., when $\la$ is not a weight of ${\rm GL}(F)$. But
in that case $K_\la(F) = 0$. While speaking of Weyl modules, we will
always tacitly assume that the rank of $F$ is big enough so that the
involved partitions are weights for ${\rm GL}(F)$ and hence the
Weyl modules in question are nonzero. With this proviso, the rank of
$F$ will not really matter to us. See the discussion after Proposition
1.1 below.

\vskip .3cm

Let us introduce the ``skew Weyl modules" defined in [ABW]. For an
arbitrary {\it skew partition} ${\la/\mu}$, [ABW] defines a functor
$K_{\la/\mu}$ assigning to any finitely generated free abelian group
$F$ the skew Weyl module $K_{\la/\mu}(F)$. (Again we will usually drop
the $F$ in the notation.) Here $\la$ is a partition and $\mu$ is any
partition whose Young diagram is contained in that of $\la$. The skew
partition ${\la/\mu}$ is best visualized as the diagram obtained by
removing the diagram of $\mu$ from the diagram of $\la$. Thus any
partition $\la$ is also a skew partition $\la/\mu$ with $\mu$ empty
and in fact ordinary Weyl modules are special cases of skew Weyl
modules.

\vskip .3cm

Given the importance of skew Weyl modules in our proofs, let us describe
briefly the nature of their definition and the fundamental results about
them in [ABW]. Only a bare sketch is given here. A thorough discussion
can be found in [ABW, section II] to which we will refer freely.
Let $D_t$ and $\Wedge_t$ respectively denote the $t$-fold divided and
exterior power functors. So for example, $D_t (F)$ is the $t$-fold
divided power of the defining representation $F$ of ${\rm GL}(F)$. For a skew
partition $\la / \mu$, let
$$D_{\la / \mu} = D_{\la_1 - \mu_1} \otimes D_{\la_2 - \mu_2} \otimes \ldots$$
and similarly for exterior powers. Then [ABW] defines $K_{\la / \mu}$ as
the image of a generalized symmetrizer map
$$d_{\la/\mu}^\prime : D_{\la / \mu} \longrightarrow \Wedge_{\tl / \tm} .$$
A concise (though informal and somewhat imprecise) way to think of the
map $d_{\la/\mu}^\prime$ is that it is ``comultiplication in the divided power
algebra along the rows of $\la / \mu$ followed by multiplication in
the exterior algebra along the columns of $\la / \mu$."

\vskip .3cm

[ABW, Theorem II.3.16] proves two important results that permit us to
think about the skew Weyl modules in two different ways. First, this
theorem explicitly describes a ``standard basis" $\{ d_{\la/\mu}^\prime (X_T) \}$
for $K_{\la / \mu}(F)$ in terms of standard tableaux $T$ of shape
$\la / \mu$. (For us, entries in a standard tableau---taken from a basis
of the defining representation $F$---will increase weakly along rows
and strictly along columns. Note that [ABW] calls such tableaux
``co-standard.") The second important result is a description of
$K_{\la / \mu}$ by generators and relations in $D_{\la / \mu}$. See
[ABW, pp. 234-5 and pp. 226-9] for details. Both these results will be
crucial in proving the main results in this paper.

\vskip .3cm

Now let us record an important general fact about the existence of
certain filtrations of skew Weyl modules proved independently in
[Kouwenhoven] and [Boffi].

\vskip .3cm

F{\sevenrm ILTRATION} T{\sevenrm HEOREM} {\it Any skew Weyl module $K_{\la / \mu}$
has a characteristic-free filtration such that the filtration factors
are isomorphic to ordinary Weyl modules, i.e., those corresponding
to partitions.}

\vskip .3cm

Such a filtration is called a Weyl filtration. Note that tensor products
of Weyl modules can also be regarded as skew Weyl modules, so the theorem
applies to them too. In the proofs of the main results we will only need
some special cases of the Filtration Theorem, in which the necessary
filtrations are constructed explicitly in [AB1, Section 3].

\vskip .3cm

{\it Stability of $K_{\la / \mu}(F)$ under variation of rank of $F$.}
The [ABW] construction of the modules $K_{\la / \mu}(F)$ ``does not
depend" on the rank $n$ of $F$. To make this precise let $E$ be a free
abelian group of finite rank $N > n$ and recall the functor $d_{N,n}$ of
[Green, Section 6.5]. $d_{N,n}$ not only takes $K_{\la / \mu}(E)$ to
$K_{\la / \mu}(F)$ but it also takes the entire apparatus involved in the
definition of $K_{\la / \mu}(E)$ (i.e., the modules and the maps) to that
involved in defining $K_{\la / \mu}(F)$. Further, $d_{N,n}$ takes a Weyl
filtration of $K_{\la / \mu}(E)$ to one of $K_{\la / \mu}(F)$. (Note that
a filtration factor $K_\nu(E)$ will become 0 under $d_{N,n}$ if $n$ is
smaller than the number of rows in $\nu$.) Perhaps a better way to express
the ``irrelevance of $F$"
is to note the following. All the constructions in [ABW] and [Boffi] are
functorial in $F$, i.e., all the involved objects are functors and all the
involved maps are natural transformations between functors. So these
constructions are the ``same" regardless of which $F$ (free abelian of
finite rank) they are applied to.

\vskip .3cm

We will occasionally need to mention linear duals of Weyl modules (``dual Weyl
modules" for short). The standard notation for the dual Weyl module with
largest weight $\la$ is Ind$_B^G \la$. The reason is that this module is also
obtained by first extending the character $\la$ to a Borel subgroup $B$, and
then inducing the character from $B$ to the whole group $G$, which for us is
${\rm GL}_n$. See, e.g., [Jantzen]. (Note that Jantzen uses the short notation
$H^0(\la)$, motivated by yet another construction of a dual Weyl module as the
space of sections of a line bundle on the flag variety.) Again we will deviate from
the standard notation and follow [ABW] by using $L_\tl(F)$, or just $L_\tl$,
to denote Ind$_B^G \la$. Note that [ABW] calls these Schur modules. A discussion
parallel to the one above for Weyl modules is valid for these Schur modules
(functoriality, definition and filtrations of skew Schur modules, stability
under variation of rank of $F$, etc.). We will only need to make limited use of
ordinary dual Weyl modules. One fact we will use is that $K_\la(F)$ and $L_\tl(F)$
are contravariant duals of each other. For example, combine [ABW, Proposition II.4.1],
with an isomorphism $F \simeq F^*$.

\vskip .3cm

{\it Notes.} (1) Given a representation $V$ of a groups $G$, the linear dual $V^*$
also becomes a left $G$-module, called the contravariant dual of $V$, by
composing the natural (right) action of $G$ on $V^*$ with an antiautomorphism
of $G$. Usually this antiautomorphism is taken to be the group inverse, but we
will take it to be the transpose of a matrix in $G = {\rm GL}_n$. This choice ensures
that the contravariant dual of any polynomial representation of a certain degree is
also a polynomial representation of the same degree. The notion of a polynomial
representation is recalled below. (2) Taking the transpose of a matrix requires choosing
an isomorphism $F \simeq F^*$, as we did above, so the contravariant duality we
use is not functorial in $F$, unlike the formulation in [ABW, Proposition II.4.1]

\vskip .3cm

Let us now discuss the Schur algebra $S(n,r)$ and its connection with the
representation theory of ${\rm GL}_n$ = ${\rm GL}(F)$. A representation $V$ of ${\rm GL}_n$
is called a polynomial representation (of degree $r$) if the entries of the
matrix by which an arbitrary $g \in {\rm GL}_n$ acts on $V$ are polynomial
functions (of degree $r$) of the entries of the matrix $g$. $K_{\la/\mu}(F)$ is
a polynomial representation of ${\rm GL}(F)$ of degree $|\la|-|\mu|$ (similarly for
Schur modules). The full subcategory of polynomial representations of ${\rm GL}_n$
of degree $r$ is equivalent to the category of left modules of the Schur algebra
$S(n,r)$. See [Green, Chapter 2] for a discussion. So Weyl and Schur modules can
be regarded as modules over Schur algebras $S(n,r)$ for appropriate $r$. The
main use of this equivalence for us will be in analyzing Ext groups, which will
be our next topic.

\vskip .3cm

{\bf Ext groups.} Our primary interest will be in the groups
Ext$^{i}_{{\rm GL}_n({\bf Z})}(K_{\la}, K_{\mu})$ for dominant weights $\la$ and
$\mu$. In the course of arguments, we will more generally need
to consider groups of the type Ext$^{i}_{{\rm GL}_n({\bf Z})}(K_{1}, K_{2})$, where
$K_1$ and $K_2$ are representations having Weyl filtrations. When no confusion
is likely, the subscript ${\rm GL}_n({\bf Z})$ will be dropped from the notation.

\vskip .3cm

Here is how the Schur algebras $S(n,r)$ enter the picture. While computing the
groups Ext$^{i}(K_{\la}, K_{\mu})$, we can (and will in all the main proofs)
assume the following without loss of generality. We can take $\la$ and $\mu$ to
be partitions (by tensoring with a suitable power of the determinant) with the
same number of boxes (since the action of the center of ${\rm GL}_n$ breaks up the
category of ${\rm GL}_n$-modules into a direct sum by degree). For the more general
Ext$^{i}(K_{1}, K_{2})$ that we need to consider, the same reasoning allows us
to assume that $K_1$ and $K_2$ are polynomial representations of the same
degree $r$. Since such $K_1$ and $K_2$ can be regarded as modules over the
algebra $S(n,r)$, we may also contemplate the group
Ext$^{i}_{S(n,r)}(K_{1}, K_{2})$. Happily [Donkin1] proves that the two kinds
of Ext groups coincide. In fact Donkin proves this result for any Schur algebra,
not just $S(n,r)$, and for any representations $K_1$ and $K_2$ of that Schur
algebra. Note further that even though [Donkin1] proves the result over algebraically
closed fields, the result extends easily to the integral setting by universal
coefficients, as made explicit in [Kulkarni1]. We will use Donkin's theorem without
further comment. In particular note that any argument about Ext groups for ${\rm GL}_n$
that requires considering the involved representations as modules over a Schur
algebra will need to use Donkin's theorem. Let us now gather several results of a
general nature about the Ext groups of interest.

\vskip .3cm

To start with let us remark that for finitely generated representations $M$ and
$N$, the groups Ext$^i(M,N)$ for $ i>0$ are all finite. Proving this will involve
a comparison with extension of scalars to {\bf Q}. The stated result follows
from Donkin's theorem (which shows that the Ext groups are finitely generated,
e.g., on account of the appropriate Schur algebra being free of finite rank over
{\bf Z}) combined with semisimplicity of representations over {\bf Q} and the
Universal Coefficient Theorem. For a statement of this theorem for group schemes
see [Jantzen, I.4.18a] and for a Schur algebra formulation, see [AB2, Theorem 5.3].
(Applying the Universal Coefficient Theorem requires $M$ and $N$ to be {\bf Z}-free,
but one easily reduces to this case by taking the quotient of $M$ and $N$ by their
respective torsion submodules $M_{tor}$ and $N_{tor}$, and then considering the
appropriate long exact sequences of Ext groups.) We will also need to use the
Universal Coefficient Theorem in some other places to relate integral Ext groups
with modular ones.

\vskip .3cm

Let us now record an important result due to Cline-Parshall-Scott-van der Kallen
[CPSvdK, Corollaries 3.2 and 3.3].

\vskip .3cm

V{\sevenrm ANISHING} T{\sevenrm HEOREM}. (Cline-Parshall-Scott-van der Kallen)
{\it $\Ext^i(K_\la, K_\mu) = 0$ unless $\la < \mu$ or $(\la = \mu$ and $i=0)$.}

\vskip .3cm

{\it Remarks.} (1) The importance of this result for the knowledege of Ext groups is
clear. In this paper we will need to use the result only in the case $\la = \mu$.
(2) The integral result as formulated above is not stated explicitly in [CPSvdK], but
it is immediate from the modular result stated and proved there. Alternatively, by
imitating later arguments of Cline-Parshall-Scott or Irving, one can give an easy
direct proof of the result over {\bf Z} using the Schur algebra setting. These authors
prove similar vanishing results (over fields) respectively in the general setting of
``highest weight categories" [CPS] and ``BGG categories" [Irving]. See
[CPS, Lemma 3.8b] and [Irving, Proposition 4.4].

\vskip .3cm

A key technical tool used to handle Ext groups in this paper will be the
the following theorem from [Kulkarni1].

\vskip 0.3 cm

S{\sevenrm KEW} R{\sevenrm EPRESENTATIVE} T{\sevenrm HEOREM.} {\it For any polynomial
representation $X$ and for a partition $\la$ containing a partition $\mu$,}
$$ {\rm Ext}^i (K_\la,L_{\tm} \otimes X) \, \simeq \, {\rm Ext}^i
(K_{\la / \mu}, X) .$$

{\it Remark.} Note the necessity of the tacit assumption on $n$ mentioned earlier,
namely $n \geq$ number of rows in $\la$. If this assumption were false,
$K_\la$ would vanish, whereas the right hand side need not. This is because
$\la / \mu$ may well have fewer rows than $\la$, possibly leading to
nonzero Ext groups on the right hand side. (The necessary assumption
on $n$ is made for the whole paper at the end of the second paragraph
on p. 253 in [Kulkarni1], but it should have been included in the statements
of the theorems there as well.) The ``irrelevance of $n$" for the results
in this paper is discussed after the proof of the next proposition.

\vskip .3cm

The next proposition is a slight strengthening of a special case of the
following result. [Green, Section 6.5] proves that if $N > n \geq r$, then
the functor $d_{N,n}$ is an equivalence of categories of modules over the
Schur algebras $S(N,r)$ and $S(n,r)$. In particular it preserves Ext groups
over these Schur algebras and by Donkin's theorem also over the corresponding
general linear groups.

\vskip .3cm

P{\sevenrm ROPOSITION} 1.1. {\it For partitions $\la$ and $\mu$,
let $N > n \geq$ the number of rows in $\la$ and number of rows in $\mu$.
Let $E$ and $F$ be free abelian groups of rank $N$ and $n$ respectively.
Then we get the following isomorphisms via the functor $d_{N,n}$ of
[Green, Section 6.5].}
$$\Ext^i_{{\rm GL}(E)}(K_\la(E), K_\mu(E)) \simeq \Ext^i_{{\rm GL}(F)}(K_\la(F), K_\mu(F)) .$$

{\it Proof.} [AB2, Section 4] gives a resolution
of $K_\la(E)$. The terms of this resolution are direct sums of modules of
the form $D_\nu(E)$ (i.e., direct sums of tensor products of divided
powers of $E$), where $\nu$ has at most as many nonzero parts as $\la$.
The construction of this resolution is functorial. So the preceding
sentences in this proof remain valid after replacing $E$ by $F$ and
in fact the two resolutions are the ``same." Speaking more precisely,
the functor $d_{N,n}$ applied to the [AB2] resolution of $K_\la(E)$
gives the [AB2] resolution of $K_\la(F)$. In particular $D_\nu$ for
the same $\nu$ and the same direct sums appear in either resolution.

\vskip .3cm

Note that by the hypothesis on $N$ and $n$, all the $\nu$ that appear in
the terms of either resolution are weights for ${\rm GL}(E)$ as well as for ${\rm GL}(F)$.
(In fact we can arrange all $\nu$ to be $\geq \la$ under the dominance
order, but we don't need this stronger fact.) This ensures that all the
$D_\nu$ are projective modules over the appropriate Schur algebra,
see [AB2]. Thus we may compute the Ext groups under consideration by
applying the appropriate Hom to these two resolutions. Now the functor
$d_{N,n}$ gives a chain map from the Hom complex involving $E$ to that
involving $F$. We can see as follows that this map is an isomorphism.
First, by [AB2, Section 2], Hom$(D_\nu(F), K_\mu(F))$ is isomorphic to
the $\nu$-weight submodule of $K_\mu(F)$; in fact this isomorphism is
natural so in particular the statement holds if we replace every $F$ by $E$.
Secondly the functor $d_{N,n}$ does not change a $\nu$-weight submodule
if $\nu$ is a weight for ${\rm GL}_n$ (and annihilates it otherwise, but this
does not happen in our situation). This completes the proof of Proposition
1.1.

\vskip .3cm

{\it Remark.} The above result is true by the same argument if one
replaces $K_\la$ and $K_\mu$ by skew Weyl modules (or even more generally,
functors with finite Weyl filtrations). $n$ should be taken at least
as much as the number of nonezero rows in the corresponding skew partitions
(for the more general version, at least as big as the number of rows
in the partitions corresponding to all Weyl modules that occur as
filtration factors in the two modules). One just has to construct suitable
projective resolutions for these modules from such resolutions for the
individual filtration factors. (Or, in case of skew Weyl modules, use the
projective resolutions already constructed in [AB2].)

\vskip .3cm

Beyond its use in proving Proposition 1.2, Proposition 1.1 is not really
necessary for the rest of the paper. But let us use the logical opportunity
it provides to discuss the minimal relevance of the $n$ in ${\rm GL}_n$ for the
main results in this paper regarding Ext groups between Weyl modules.
Proposition 1.1 is not enough by itself to address this issue because the
proofs will more generally employ Ext groups involving modules with Weyl
filtrations. In view of Green's result quoted before Proposition 1.1, we
can a priori take $n$ to be large enough. After a ``stable" result about Ext
groups between two Weyl modules is proved, Proposition 1.1 then extends such
a result to all $n$ for which the partitions involved in the statement of the
result are weights for ${\rm GL}_n$. But a better reasoning is to simply note that
the proofs of the main results stay valid for any such ``appropriate" value
of $n$. (For details on how various ingredients involved in the proofs behave
vis a vis the value of $n$, see the Remark after Proposition 1.1, the Remark
after the statement of the Skew Representative Theorem above, and the paragraph
discussing stability of skew Weyl modules under variation of rank of $F$.)

\vskip .3cm

Before turning to Proposition 1.2 we need to record the following
isomorphisms of Ext groups. First, by contravariant duality,
$\Ext^i(K_\la,K_\mu) \simeq \Ext^i(L_{\tm}, L_{\tl}).$
Combining this with the Ext-preserving functor $\Omega$ (a weak form of
Howe duality) from [AB2, Section 7], one gets the following conjugate
symmetry of Ext groups. Ext$^i(K_\la, K_\mu) \simeq \Ext^i(K_{\tm}, K_{\tl})$.
Note that one may need to increase $n$ (valid by Proposition 1.1) so
that the conjugate partitions are weights for ${\rm GL}_n$.

\vskip .3cm

P{\sevenrm ROPOSITION} 1.2. (column and row removal principles) {\it Suppose
the lengths of the first columns (rows) of partitions $\la$ and $\mu$
are equal. Let $\la ^\prime$ and $\mu^\prime$ be the partitions obtained
by removing the first columns (rows) of $\la$ and $\mu$ respectively. Then}
$\Ext^i(K_\la, K_\mu) \simeq \Ext^i(K_{\la ^\prime}, K_{\mu ^\prime})$.

\vskip .3cm

{\it Proof.} Let us first prove the column removal principle.
Call the common length of the first columns $m$. Evidently $m \leq n$.
By Proposition 1.1 we can change $n$ to $m$. Then
$K_\la = K_{\la ^\prime} \otimes det, K_\mu = K_{\mu ^\prime} \otimes det$,
where {\it det} is the one-dimensional determinant representation.
The claim is immediate by canceling the determinant. Row removal
now follows by conjugate symmetry of Ext groups. (It may be necessary
to increase and then restore $n$, again using Proposition 1.1.)

\vskip .3cm

{\it Remark.} Results similar to the two Propositions above, but about
decomposition numbers in characteristic $p$, have been proved in
[Green] and [James] respectively. Let us outline an approach to these
earlier results using certain Ext groups, namely those between
Weyl modules and simple modules in characteristic $p$. The analogues
of Proposition 1.1 and the column removal principle for these ``mixed"
Ext groups are easily obtained by repeating the arguments above.
Combining this with the well-known connection between Euler
characteristics built out of such Ext groups and decomposition numbers,
one gets [Green, Theorem 6.6e] and [James, Theorem 1]. For [James, Theorem
2] (row removal principle for decomposition numbers), additional arguments
seem necessary, as the functor $\Omega$ is not available in this situation.
For instance one can use the connection with tilting multiplicities
obtained by Donkin to get the desired result. Let us skip the details.

\vskip 1cm

\centerline{2. ON THE EXTENSIONS BETWEEN WEYL MODULES}

\vskip .5cm

This section contains the main results of this paper concerning
Ext groups between certain Weyl modules for ${\rm GL}_n$ over the integers.
The significance of these results for modular representation theory is
discussed in the Remarks after each result.

\vskip .3cm

Let us outline the general strategy that we will use to study
Ext groups between Weyl modules. As remarked before, one can
without loss of generality take the corresponding two dominant
weights to be partitions with the same number of boxes. Now the Skew
Representative Theorem [Kulkarni1] gives the following isomorphism
$$ {\rm Ext}^i (K_\la,L_{\tm} \otimes K_\nu) \, \simeq \, {\rm Ext}^i
(K_{\la / \mu}, K_\nu) , $$
where $\la, \mu, \nu$ are partitions with the diagram of $\mu$ contained
in that of $\la$ and $|\nu| = |\la|-|\mu|$. (Let me again mention that
following [ABW], $L_{\tm}$ denotes the dual Weyl module with highest weight
$\mu$.) Notice that the degree of representations occuring on the LHS is
$|\la|$ whereas that on the RHS is $|\la|-|\mu|$. This leads to the hope
of somehow analyzing Ext groups recursively, using the just mentioned
reduction in degree. A difficulty involved in this approach is that one
has to deal with more complicated representations, namely $K_{\la / \mu}$
and $L_{\tm} \otimes K_\nu$. One can try to deal with  these two objects
piece by piece by working with their filtrations. One does have a Weyl
filtration of $K_{\la / \mu}$ which, however, is very complicated in
general.  Worse still, no information in general is known about the
filtrations of the latter.

\vskip .3cm

One can get around both these difficulties by choosing $L_{\tm}$ to be
an exterior power $\Wedge_t$ of the defining representation, i.e., by
taking $\mu = 1^t$, a single column of length $t$. Since $\Wedge_t$ is
a Weyl module as well, $\Wedge_t \otimes K_\nu$ also has a Weyl
filtration. Weyl filtrations of both these modules (i.e., of $K_{\la / 1^t}$
and $\Wedge_t \otimes K_\nu$) are constructed explicitly in [AB1, Section 3].
These filtrations are fairly simple to describe and one has a very easy
description of the Weyl modules that occur as filtration factors by
Pieri-type rules. In the simplest new case (treated in Theorem 2.1) one
can inductively control the maps induced on Ext groups as we patch together
the Weyl modules in these filtrations. In order to seriously entertain the
idea of computing {\it all} Ext groups between Weyl modules using this
technique, one must be able to grapple successfully with spectral sequences
arising from Ext groups applied to filtered complexes whose terms have Weyl
filtrations. This seems difficult in general. Here is what one can do at
present using the approach outlined above.

\vskip .3cm

(1) Theorem 2.1 gives a complete determination of the groups
Ext$^{i}(K_{\la}, K_{\mu})$ where ${\mu}-{\la}$ is a positive root $\alpha$ of
${\rm GL}_n$. This is the simplest case left open after the [CPSvdK] vanishing theorem.
(2) For general $\la < \mu$, one can circumvent the difficulty involved in
dealing with filtrations that was mentioned above by attempting a less
ambititious task. One can compute a ``multiplicative Euler characteristic"
$\chi$ defined as the alternating product of cardinalities of Ext groups
between a given pair of Weyl modules. One first gets a recursive algorithm
to do this computation which in turn leads to a simple formula for $\chi$.
This is done in [Kulkarni2]. The base of the recursion is precisely the
case treated by Theorem 2.1.
(3) The Skew Representative Theorem takes a particularly simple form if
the first partition is a single column. Using this, one can show that
Ext$^1(\Wedge_{|\la|}, K_\la)$ is always cyclic. A little more reasoning
enables one to explicitly calculate the order of this Ext$^1$ in terms of
$\la$. By contravariant duality and conjugate symmetry of Ext groups, one
also determines three other types of $\Ext ^1$ groups. This is the content
of Theorem 2.2.

\vskip .3cm

Let us tackle Theorem 2.2 first as its proof is much shorter and gives a quick
demonstration of the use of two tools that will be used much more substantially
while proving Theorem 2.1: (1) the Skew Representative Theorem and
(2) the [ABW] description of Weyl modules by generators and relations.
Theorem 2.2 computes the group Ext$^1(\Wedge_{|\la|}, K_\la)$. Note that if the
partition $\la$ has only one column, then $K_\la = \Wedge_{|\la|}$ and
this Ext group is trivial by the [CPSvdK] vanishing theorem. So let us
assume that $\la$ has at least two columns.

\vskip .3cm

T{\sevenrm HEOREM} 2.2. {\it The group $\Ext ^1 (\Wedge_{|\la|}, K_\la)$
between an exterior power of the defining representation and another
Weyl module is always cyclic and its order is the {\it gcd g} of integers
$${a+1 \over gcd(a+1, lcm(1,2,3,\ldots,b))} \; ,$$
where $a,b$ are the lengths of any two consecutive columns of $\la$ with
$a \geq b$.}

\vskip .3cm

C{\sevenrm OROLLARY.} {\it The statement of Theorem 2.2 stays valid if we replace
the Ext group there by any of the following: $\Ext ^1 (L_\tl, \Wedge_{|\la|}),
\, \Ext ^1 (K_\tl, D_{|\la|})$ and \/ $\Ext ^1 (S_{|\la|}, L_\la)$,
where $D$ and $S$ denote respectively the divided and symmetric powers of
the defining representation.} (This follows from the fact that the four
Ext groups in question are isomorphic by contravariant duality and
conjugate symmetry.)

\vskip .3cm

{\it Proof of Theorem 2.2.}
If $\la$ has exactly two columns of lengths $a \geq b$, then by [AB2, Section 9]
$\Ext ^1 (\Wedge_{|\la|}, K_\la)$ is cyclic and its order is as given above.
If $\la$ has more than two columns, pick one of them, say of length $\ell$,
and call the partition left by erasing that column $\la '$. Take the exact
sequence $0 \arr K_\la \arr K_{\la '} \otimes \Wedge_\ell \arr X \arr 0$
obtained from the Weyl filtration of $K_{\la '} \otimes \Wedge_\ell$ given
in [AB1, Section 3]. Apply Hom$(\Wedge_{|\la|},-)$ and take the associated
long exact sequence. It begins as follows. $0 \arr \Ext ^1 (\Wedge_{|\la|}, K_\la)
\arr \Ext ^1 (\Wedge_{|\la|}, K_{\la '} \otimes \Wedge_\ell)$. By the Skew
Representative Theorem, $\Ext ^1 (\Wedge_{|\la|}, K_{\la '} \otimes \Wedge_\ell)
\simeq \Ext ^1 (\Wedge_{|\la '|}, K_{\la '})$. Now we can repeat the process
until only two columns are left. Clearly these can be arranged to be any two
columns of the original partition $\la$. So we have proved that the desired
Ext$^1$ is cyclic and that its order divides $g$. To show that the order
equals $g$, we will instead use Ext$^1 (K_\tl, D_{|\la|})$, which is isomorphic
to the desired Ext$^1$.

\vskip .3cm

Let us recall some facts from [ABW, Section II] about $K_\tl$ (see also Section 3,
specifically the beginning of the proof of Lemma C and Note 2 at the end of Example
2 in the same proof). We have $d_\tl ': D_\tl \arr \hskip -.2cm \arr K_\tl$. Label a
basis of the defining representation as $e_1, e_2, \ldots$ . Then we have a cyclic
generator $e_1^{(\tl_1)} \otimes e_2^{(\tl_2)} \otimes \ldots$ for $D_\tl$ (and
hence one for $K_\tl$) and the kernel of $d_\tl'$ is generated by elements
$$\ldots \otimes e_k^{(a)} \otimes e_k^{(t)} e_{k+1}^{(b-t)} \otimes \ldots ,
\quad 1 \leq t \leq b ,$$
where $k, k+1$ are any two consecutive rows of the partition $\tl$ of lengths $a \geq b$
respectively and ``$\ldots$" indicates tensor factors $e_i^{(\tl_i)}$ in all positions
$i$ other than $k$ and $k+1$.

\vskip .3cm

By the Universal Coefficient Theorem it suffices to show that $g$ is the largest
integer modulo which there is a nonzero equivariant map $K_\tl \arr D_{|\la|}$.
So let us work over ${\bf Z}/m{\bf Z}$ and characterize the integers $m$ for which
such a map exists. Since all weight spaces of $D_{|\la|}$ are one-dimensional,
Hom$(D_\tl, D_{|\la|})$ is generated by the map taking the cyclic generator
$e_1^{(\tl_1)} \otimes e_2^{(\tl_2)} \otimes \ldots$ to the basis element
$e_1^{(\tl_1)} e_2^{(\tl_2)} \ldots$ of the relevant weight space. Clearly this
map is just multiplication in the divided power algebra. This map will descend to
$K_\tl$ exactly when it kills the generators of the kernel of $d_\tl'$ listed above.
Recalling how to multiply in the divided power algebra, this in turn will happen
exactly when ${a+t \choose t} = 0$ in ${\bf Z}/m{\bf Z}$ for all $1 \leq t \leq b$,
where $a \geq b$ are lengths of any two consecutive rows of $\tl$, i.e., consecutive
columns of $\la$. Fixing $a$ and $b$, the $gcd$ of the resulting $b$ binomial
coefficients is easily seen to be the number displayed in the statement of Theorem
2.2, thus completing the proof.

\vskip .3cm

{\it Remarks.}
(1) It should be clear that the argument in the previous paragraph can prove all of
Theorem 2.2, making the first paragraph of the proof redundant. (The first paragraph
was nevertheless included to illustrate how the Skew Representative Thorem may be
used.) Such a streamlined proof would also recover as a special case the result from
[AB2, Section 9] quoted in the first paragraph. Actually the argument in the previous
paragraph is close to the one in [AB2], just formulated differently and applied in
a more general situation.

\vskip .3cm

(2) In fact one can get an even more general result from the argument in the last
paragraph of the proof. Since {\it skew} Weyl modules are also cyclic and have a
description by generators and relations, one can easily calculate by the same
reasoning Ext$^1$ between any skew Weyl module and a divided power of the defining
representation (and three other types of Ext groups by the argument in the Corollary).
This group is also cyclic and its order is equal to the {\it gcd} of several binomial
coefficients determined by lengths of and overlap between successive pairs of adjacent
rows of the skew partition in question. (The simplification contained in the last
sentence of the proof will generally not be available here.)

\vskip .3cm

(3) Let us work in prime characteristic $p$ in this and the next paragraph.
Again using the Universal Coefficient Theorem, the Corollary to Theorem 2.2
shows that Hom$(S_{|\la|}, L_\la)$ always has dimension 0 or 1, and that it
is one-dimensional exactly when $p$ divides the number appearing in Theorem
2.2 for each pair of successive column lengths $a,b$ of the partition $\la$.
Thus we get a necessary and sufficient condition for the existence of a
homomorphism from a symmetric power of the defining representation to a dual
Weyl module $L_\la$ of largest weight $\tl$. This condition is easier to test
by writing the length of each row of weight $\tl$ (we are expressing everything
in terms of the highest weight $\tl$) $p$-adically and displaying the digits for
successive rows of $\tl$ in a kind of tableau form. Thus the $i$-th column from
right of this ``digit tableau" displays the $p^i$-place digits of successive row
lengths of $\tl$. The condition obtained from Theorem 2.2 then amounts to the
requirement that in the digit tableau each digit strictly above and weakly to the
right of any nonzero digit must be $p-1$.

\vskip .3cm

Note that the existence of a homomorphism from a symmetric power of the defining
representation to a dual Weyl module $L_\la$ implies in particular that the
simple module with highest weight $\tl$ is a composition factor of that symmetric
power, though of course the former is a significantly stronger requirement. It is
interesting to compare the composition factors of $S_{|\la|}$ -- known as a result
of the known ${\rm GL}_n$-submodule structure of symmetric powers [Doty] -- with the
condition in Theorem 2.2. This is easier to state using the $p$-adic digit tableau
above. (This language is modelled on that in [Krop], which also proves the same
submodule structure for symmetric powers under the action of the full matrix
semigroup. The idea of exploiting one-dimesionality of weight spaces in symmetric/divided
powers too was taken from the independent works of Doty and Krop.) The highest weights
of the composition factors of $S_{|\la|}$ are characterized by each column of digits
in the digit tableau of the highest weight having the form of a string of $(p-1)$'s
possibly ending with a single different digit. Clearly this is a weaker condition
than the one in the previous paragraph. (Note that we always take the $n$ in
${\rm GL}_n$ as large as it needs to be, so the condition on the size $n$ in the
enumeration of composition factors of $S_{|\la|}$ is automatically met.)

\vskip .3cm

T{\sevenrm HEOREM} 2.1. {\it Suppose $\la$ and $\mu$ are dominant weights for ${\rm GL}_n$
such that ${\mu}-{\la}$ is a positive root $\alpha$. Let $\rho$ be half
the sum of positive roots. Then $\Ext^{i}(K_{\la}, K_{\mu})$ is
cyclic of order ${\langle} {\la} + {\rho} , {\alpha} \, \check{} \, {\rangle} + 1$
for $i = 1$ and vanishes for all other $i$.}

\vskip .3cm

To facilitate the proof of Theorem 2.1, let us translate this statement into
combinatorial language and make some reductions. We have $\la = \sum {\la_i \epsilon_i}$,
where $\epsilon_i$ is the weight consisting of 1 in the $i$-the position and
zeroes elsewhere. Let $\alpha$ be the positive root $\epsilon_r - \epsilon_s$ for
$r<s$. Then the claimed order of Ext$^{1}(K_{\la}, K_{\mu})$ is $\la_r - \la_s + s - r + 1$.
As noted before, without loss of generality we can take $\la$ and $\mu$ to be
partitions with the same number of boxes. Then the Young diagram of $\mu$ is
obtained from that of $\la$ by raising a single box from the end of the $s$-th
row up to the end of the $r$-th row. The number $\la_r - \la_s + s - r + 1$ is easily
seen to be the hook length of the box in the $r$-th row and $\la_s$-th column of
the diagram of $\la$ (or the diagram of $\mu$). (Recall that the hook length of a
box in the Young diagram of a partition is the total number of boxes to the right
and below it, including itself.) Now by appealing to the row and column removal
principles, we may further strip off the $r-1$ identical first rows and $\la_s-1$
identical first columns from $\la$ and $\mu$ without changing the Ext groups.
Then the assertion about Ext$^{1}(K_{\la}, K_{\mu})$ amounts to saying that this
group is cyclic of order equal to the hook length of the box in the top left corner
of the diagram of $\la$ (or $\mu$). Let us set up some notation so we can
precisely state and prove Theorem 2.1 in this equivalent combinatorial form.

\vskip .3cm

Let $a_i, p_i$ for $i= 1,\ldots, k$ be positive integers with
$a_1 > a_2 > \ldots > a_k $. Let $\la = a_1^{p_1} \ldots a_k^{p_k}1$,
$\mu = (a_1+1) a_1^{p_1-1} a_2^{p_2} \ldots a_k^{p_k}$ and
$\nu =  a_1^{p_1} \ldots a_k^{p_k}$.
In words, $\nu$ is a partition whose diagram consists of $k$ rectangular
blocks of rows; $\la$ and $\mu$ are obtained by adding a single box to the
first column and the first row of $\nu$ respectively. Let $h_i$ be the hook
length of the $(p_1+\ldots+p_{i-1}+1)$-th box from top in the first column
of $\la$, i.e., the top left box in the $i$-th rectangular block in the
diagram of $\la$. Let $\ell_i$ be the hook length of the $(p_1+\ldots+p_i)$-th
box from top in the first column of $\nu$, i.e., the bottom left box in
the $i$-th rectangular block in the diagram of $\nu$. For future use, note
that $\ell_j = h_j - p_j =  a_j + p_{j+1} + p_{j+2} + \ldots + p_k$. Now we
can state and prove an equivalent form of Theorem 2.1 using this notation.
For technical convenience we will exclude the Hom case from the statement.
(See the paranthetical note at the end of the first paragraph of the proof.)
We may do so since it is already known that ${\rm Hom}(K_\la, K_\mu) = 0$
for $\la \neq \mu$.

\vskip .3 cm

T{\sevenrm HEOREM} 2.1. (Combinatorial version) {\it For $\la, \mu$ as in the
preceding paragraph,}
$$\Ext^1(K_\la, K_\mu) = {\bf Z}/h_1{\bf Z} , \quad
\Ext^i(K_\la, K_\mu) = 0 \hbox { for } i > 1 .$$

{\it Proof.} Note that $\mu$ is obtained by removing the
single box from the last row of $\la$ and placing it at the end
of the first row. By induction we assume the result for all such
pairs $(\la , \mu )$ of smaller degree. This induction starts in
degree 1. Here $\nu$ must be empty, $\la = \mu$ = a single box and
$h_1 = 1$. We need to show that all Ext$^i$ vanish for $i>0$, which
is immediate by the [CPSvdK] vanishing theorem. (This case does not
exist in the original formulation of Theorem 2.1, but does make sense
for the combinatorial version we are proving. Also note that in this
case Hom is nonzero, which is why the Hom case was excluded in the
statement of the theorem. Alternatively we could have started the
induction in degree 2, where again there is just one case, namely
when $\nu$ consists of a single box. Then one needs to show that
$Ext^i(\Wedge_2, D_2)$ is ${\bf Z}/2{\bf Z}$ for $i=1$ and $0$
otherwise. This is immediate, e.g., from the projective resolution
$0 \arr D_2 \arr F \otimes F \arr \Wedge_2 \arr 0$.
But in fact this case is subsumed in the inductive step below.)

\vskip .3cm

We will use the explicit Weyl filtrations for $K_{\la/1}$ and
$K_\nu \otimes F$ constructed in [AB1, Section 3]. The filtration
factors are described by Pieri-type rules. The factors for $K_{\la/1}$
correspond to partitions obtained by deleting one box from the diagram
of $\la$. The factors for $K_\nu \otimes F$ correspond to partitions
obtained by adding one box to the diagram of $\nu$.

\vskip .3cm

The Skew Representative Theorem implies that 
$\Ext^i(K_\la,  F \otimes K_\nu) \simeq \Ext^i(K_{\la/1}, K_\nu)$. So let
us consider the [AB1] Weyl filtrations of $K_{\la/1}$ and $F \otimes K_\nu$.
These give us the following exact sequences:
$0 \arr K_\nu \arr K_{\la/1} \arr M \arr 0$,
$0 \arr K_\mu \arr F \otimes K_\nu \arr N \arr 0$ and
$0 \arr P \arr N \arr K_\la \arr 0$,
where the Weyl modules occuring in $grM$, $grN$ and $grP$ are prescribed by
Pieri's rules. We will analyze certain long exact sequences associated to
these short exact sequences.

\vskip .3cm

First apply Hom$(-,K_\nu)$ to $0 \arr K_\nu \arr K_{\la/1} \arr M \arr 0$
and take the corresponding long exact sequence. By [CPSvdK],
Ext$^i(K_\nu, K_\nu) = 0$ for $i \neq 0$. Our inductive hypothesis combined
with the row removal principle applies as we patch together the Weyl
modules occuring in $grM$ and take the corresponding long exact sequences
of Ext groups with $K_\nu$. We get as a result that Ext$^i(M, K_\nu) = 0$
for $i \neq 1$ and that the cardinality of Ext$^1(M, K_\nu)$ equals
$\ell_1 \ldots \ell_k$. Known information now forces Ext$^i(K_{\la/1}, K_\nu) = 0$
for $i > 1$. Therefore the long exact sequence reduces to
$$0 \larr {\rm Hom} (K_{\la/1}, K_\nu)  \larr {\rm Hom} (K_\nu, K_\nu)
\larr \Ext^1(M, K_\nu) \larr \Ext^1 (K_{\la/1}, K_\nu) \larr 0 . $$
As $K_\nu$ occurs once in $gr (K_{\la/1})$, the first two Hom terms are both
isomorphic to {\bf Z}. Therefore the map between them is given by an integer.
Now we state

\vskip .3cm

L{\sevenrm EMMA} A. {\it The map ${\rm Hom} (K_{\la/1}, K_\nu)  \arr {\rm Hom} (K_\nu, K_\nu)$
induced by the inclusion $K_\nu \hookrightarrow K_{\la/1}$ arising from the
[AB1] Weyl filtration of $K_{\la/1}$ is given by the integer $\ell_1 \ldots \ell_k$.}

\vskip .3cm

Assuming Lemma A, $\Ext^1(M, K_{\nu})$ is forced to be isomorphic to
{\bf Z}/$\ell_1 \ldots \ell_k${\bf Z} and $\Ext^1 (K_{\la/1}, K_\nu)$
is seen to vanish. To sum up, modulo Lemma A, we have proved that
$\Ext^i (K_{\la/1}, K_\nu) = 0$ for $i \neq 0$. Therefore
$\Ext^i (K_\la, F \otimes K_\nu) = 0$ for $i \neq 0$ by the
Skew Representative Theorem. 

\vskip .3cm

Next we apply Hom$(K_\la, -)$ to $0 \arr P \arr N \arr K_\la \arr 0$ and
take the corresponding long exact sequence. Ext$^i(K_\la, K_\la) = 0$ for
$i \neq 0$ by the [CPSvdK] vanishing theorem. Again our inductive
hypothesis combined with row removal principle applies as we patch
together the Weyl modules occuring in $grP$ and take the corresponding
long exact sequences of Ext groups of $K_\la$ with the resulting modules.
We get as a result that Ext$^i(K_\la, P) = 0$ for $i \neq 1$ and that the
cardinality of Ext$^1(K_\la, P)$ equals $h_2 \ldots h_k$. Known information
now forces Ext$^i(K_\la, N) = 0$ for $i > 1$. Therefore the long exact
sequence reduces to
$$0 \larr {\rm Hom} (K_\la, N)  \larr {\rm Hom} (K_\la, K_\la)
\larr \Ext^1(K_\la, P) \larr \Ext^1 (K_\la, N) \larr 0. $$
As $K_\la$ occurs once in $grN$, the first two Hom terms are both
isomorphic to {\bf Z}. Therefore the map between them is given by
an integer. Now we state

\vskip .3cm

L{\sevenrm EMMA} B. {\it The map ${\rm Hom} (K_\la,N)  \arr {\rm Hom} (K_\la, K_\la)$
induced by surjection $N \arr \hskip -.2cm \arr K_\la$ arising from the
[AB1] Weyl filtration of $F \otimes K_\nu$ is given by the integer
$h_2 \ldots h_k$.}

\vskip .3cm

Assuming Lemma B, $\Ext^1(K_\la, P)$ is forced to be isomorphic to
{\bf Z}/$h_2 \ldots h_k${\bf Z} and $\Ext^1 (K_\la, N)$ is seen
to vanish. So modulo Lemma B, we have proved that
$\Ext^i (K_\la, N) = 0$ for $i \neq 0$.

\vskip .3cm

Now we are ready to apply Hom$(K_\la, -)$ to
$0 \arr K_\mu \arr F \otimes K_\nu \arr N \arr 0$
and analyze the resulting long exact sequence.
Clearly Hom$(K_\la, K_\mu)=0$. We already know that $\Ext^i (K_\la, N) = 0$ and
$\Ext^i (K_\la, F \otimes K_\nu) = 0$ for $i \neq 0$. Therefore
$\Ext^i(K_\la, K_\mu)$ is forced to vanish for $i>1$ and the exact sequence
reduces to
$$0 \larr {\rm Hom} (K_\la, F \otimes K_\nu)  \larr {\rm Hom} (K_\la, N)
\larr \Ext^1(K_\la, K_\mu) \larr 0. $$
As $K_\la$ occurs once as a factor in the filtrations of $F \otimes K_\nu$ as
well as $N$, the two Hom terms are both isomrphic to {\bf Z}. Therefore the
map between them is given by an integer. Now we state

\vskip .3cm

L{\sevenrm EMMA} C. {\it The map ${\rm Hom} (K_\la, F \otimes K_\nu)  \arr {\rm Hom}
(K_\la, K_\la)$ induced by the surjection $F \otimes K_\nu \arr \hskip -.2cm
\arr K_\la$ arising from the [AB1] Weyl filtration of $F \otimes K_\nu$ is given
by the integer $h_1 \ldots h_k$.}

\vskip .3cm

Now Lemma B and Lemma C together imply easily that the map
${\rm Hom} (K_\la, F \otimes K_\nu)  \arr {\rm Hom} (K_\la, N)$ is given by
the integer $h_1$. Therefore $\Ext^1(K_\la, K_\mu) \simeq {\bf Z}/h_1{\bf Z}$.
Thus Theorem 2.1 is proved modulo Lemmas A, B and C. The proofs of these
lemmas are similar to each other and somewhat intricate. They will be
presented in the next section.

\vskip .3cm

{\it Remarks.}
(1) Let us discuss the significance of Theorem 2.1 for modular representation
theory. So we will work over a field of characteristic $p$ until further notice.
Consider the dot action of the affine Weyl group $W_p$ on weights. Let $\la$ be
a dominant regular weight in an alcove $C$ and $s$ the reflection in a wall of $C$
such that $\la < s.\la$ (and so $s.\la$ is also dominant). Then $K_{\la}$ and
$K_{s.\la}$ are called neighboring Weyl modules. It is known that for any split reductive
algebraic group $\Ext ^i(K_{\la}, K_{s.{\la}})$ is one dimensional for $i = 0$ or 1
and vanishes otherwise. See, e.g., [Jantzen, II.7.19d]. The nonzero homomorphism
is of interest because, for instance, it is involved in an important conjecture of
Jantzen. See [Andersen1]. Here is how Theorem 2.1 combined with the tranlsation
principle can be used to recover this Ext calculation for ${\rm GL}_n$. (This strategy
is admittedly strange, since the stated result was obtained using just translation
functors, albeit in a little more substantial way. The point here is to observe that
the case in Theorem 2.1 is enough to realize neighboring modules between any two
given adjacent dominant alcoves.) We have $s.\la - \la =  k \alpha$,
a positive multiple of a positive root $\alpha$. Consider first the case $k = 1$,
so the two weights are as in Theorem 2.1. It is easy to see that $p$ divides
${\langle} {\la} + {\rho} , {\alpha} \, \check{} \, {\rangle} + 1$ and then
the Universal Coefficient Theorem gives the desired modular Ext groups. Now
suppose $k > 1$. Note that by the translation principle
Ext$^i (K_\la, K_{s.\la}) \simeq \Ext^i (K_\eta, K_{s.\eta})$ for any dominant
regular weight $\eta$ in $C$. We will choose $\eta$ close to the wall separating
$C$ and $s.C$ so as to reduce to the case $k = 1$. For this first note that all
weights between the two regular weights $\la$ and $s.\la$ (on the
straight line connecting them) except $\la + (k/2) \alpha$ are regular and lie
either in $C$ or $s.C$. The exception, not necessarily an integral weight,
is where the line intersects the wall separating $C$ and $s.C$. If $k = 2t + 1$
is odd, then let $\eta = \la + t \alpha$ and then $s.\eta = \la + (t+1) \alpha$,
reducing to the case $k = 1$. If $k = 2t$ is even, then we may first replace
$\la$ by $\la + (t-1) \alpha$ and so $s.\la$ by $\la + (t+1) \alpha$. In
other words we may assume that $k=2$. Let $\alpha = \epsilon_q - \epsilon_r$
with $q < r$. Now take $\eta = \la - \epsilon_r$. Then $s.\eta = s.\la - \epsilon_q$.
(If $\la$ and $s.\la$ are partitions, $k=2$ means $s.\la$ is obtained by moving two
boxes from a lower row of $\la$ to a higher row. Now $\eta$ as above is simply
obtained by erasing one of these two boxes from $\la$ and then $s.\eta$ similarly
has one box less than $s.\la$.) It is easy to check that $\eta$ is in $C$ and that
we are again in the case $k = 1$.

\vskip .3cm

A reason why Theorem 2.1 is of additional interest is as follows. The definition of
neighboring Weyl modules makes sense only if there are regular weights, which happens
iff $p \geq n$, the Coxeter number for ${\rm GL}_n$. But even for small $p$, Theorem 2.1
still gives results about Ext groups between Weyl modules whose dominant weights differ
by a single root.

\vskip .3cm

(2) Continue with neighboring Weyl modules $K_\la$ and $K_{s. \la}$ as in Remark 1
but work over ${\bf Z}_p$, the ring of $p$-adic integers. Let $v$ be the largest
power of $p$ dividing ${\langle} {\mu} + {\rho} , {\alpha} \, \check{} \, {\rangle}$,
where $\mu$ is a weight on the wall separating the alcoves containing $\la$ and $s. \la$.
About the same time when the combinatorial version of Theorem 2.1 was first proved,
Andersen proved that for any split reductive algebraic group Ext$^1(K_\la, K_{s. \la})$
is cyclic of order $p^v$ and all other Ext$^i(K_\la, K_{s. \la})$ vanish. See
[Andersen3]. For ${\rm GL}_n$, Andersen's result can be recovered from Theorem 2.1
and base change to ${\bf Z}_p$ using the same argument that was used in
Remark 1 to recover the weaker result on modular Ext. (As shown in [Andersen1]
translation functors work over ${\bf Z}_p$ too.) Again the result over
${\bf Z}_p$ cannot apply for small $p$. But Andersen's result and Theorem 2.1
together lead one to hope that the statement of Theorem 2.1 should stay valid
for any reductive algebraic group over {\bf Z}. Granting the use of translation,
such a result (if true) can be thought of as a natural common generalization of
Theorem 2.1 and Andersen's result. See the Introduction for an indication of some
possible approaches to proving such a result.

\vskip .3cm

(3) Now let $\la$ and $\mu$ be as in the statement of Theorem 2.1, but work over
a field of characteristic $p$. Theorem 2.1 contains more information than
is registered in the modular calculation of the groups Ext$^i(K_\la, K_\mu)$
indicated in Remark 1. The extra information is the power of $p$ dividing
${\langle} {\la} + {\rho} , {\alpha} \, \check{} \, {\rangle} + 1$. This power
turns out to be the coefficient of ch($K_\la$) in the Jantzen sum formula for
$K_\mu$ obtained from the determinant of the contravariant form on integral
$K_\mu$. Actually in general this coefficient is encoded in the Euler
characteristic $\chi$ defined in the discussion near the begining of this
section, but for the case in Theorem 2.1 $\chi$ amounts to Ext$^1$ due to
the vanishing of other Exts. This connection with $\chi$ is explained in
[Kulkarni2]. It is unclear if the vanishing of higher Exts in Theorem 2.1
has any direct representation theoretic significance.

\vskip .3cm

(4) The proof of Theorem 2.1 illustrates the lines along which one may hope
to calculate Ext groups in general using the Skew Representative Theorem.
In view of the discussion on neighboring Weyl modules, one case of interest
is the generalization when $\mu - \la = k \alpha$ with $k > 1$. (This is a
generalization because for neighboring Weyl modules the $k$ depends on the
weight $\la$ and the prime $p$.) The difficulties here are formidable, as
can be seen from the complexity of the answer calculated in [BF] for Ext$^1$
in this case for ${\rm GL}_3$. Unlike the simple answer given by Theorem 2.1, the
answer in [BF] involves taking the $gcd$ of many numbers (this is true even
for the ${\rm GL}_2$ case in [AB2] that was used in Theorem 2.2) and uses the lengths
of all rows of $\la$.

\vskip .3cm

In the light of quickly mounting and seemingly inherent difficulties that
one faces when $k > 1$, the fact that one has such a simple answer over {\bf Z}
for $k = 1$ (and possibly valid for other reductive groups too) acquires added
interest. It is tempting to wonder if the $k = 1$ case can at least partially
serve as a reasonable {\bf Z} analogue of the notion of neighboring Weyl
modules.

\vskip 1cm

\centerline{3. PROOFS OF LEMMAS A, B AND C}
\vskip .5cm

The strategy behind the proofs of the three lemmas is quite simple, but to
carry it out requires some notational and computational effort. Each lemma
identifies the integer giving a specific map between two Hom groups
each of which is isomorphic to {\bf Z}. The proofs explicitly find
a generator of the source Hom group (the hardest step), apply the map to this
generator, and inspect which multiple of the generator one gets in the target
Hom group. The hardest step comes down to explicitly finding and solving a
linear system of integer equations. What is noteworthy is that beyond keeping
track of all the paraphernalia required for bookkeeping, the lemmas offer no
further difficulty. The computations stay reasonable at all stages, in
particular one does not have to resort to taking {\it gcd\/}'s at any point.
This should be contrasted with what usually happens in this type of setting,
e.g., see Remark 4 after the proof of Theorem 2.1. Further, one can give
very similar (and equally long) proofs for all three lemmas and all features
occuring in such proofs (except some straightening required in one step for
Lemma A) can be already seen in the much simpler special cases of Lemma C
worked out below. Moreover, here we will be able to get away with doing the
hardest step only once (for Lemma C) and thus substantially shorten the proofs
of the remaining two lemmas.

\vskip .3cm

These facts suggest that there may be a more efficient way to organize these
computations, and perhaps even existence of a more conceptual explanation.
The following considerations might be relevant in this regard. Just like the
answers claimed in the three lemmas, several objects in the proofs are naturally
indexed by $k$-tuples, where $k$ is the number of rectangular blocks in the
partition $\nu$. It seems plausible that there should be a simpler approach to
proving these lemmas that uses induction on $k$. Such induction does make an
appearance in the proofs below, but only at the very end when one has to simplify
a laboriously obtained algebraic expression. In another direction, since the order
of the Ext$^1$ in Theorem 2.1 arises ultimately from the three lemmas, it is
interesting to look for the origin of Andersen's $p$-adic answer [Andersen3]
in his proof. In his computation, a role analogous to that of the three lemmas
is played by [Andersen2, Lemmas 2.2 and 2.4]. To prove these lemmas Andersen
gives a short and slick argument using adjointness properties of translation
functors. He then brings in the Weyl chracter formula to get the final numerical
form of the answer. It would be nice to be able to adapt his approach to the
integral situation (perhaps still working with one prime at a time, but letting
the prime to be arbitrary, i.e., allowing nonregular weights and thereby allowing
more complicated modules with Weyl filtration to arise upon translation).

\vskip .3cm

After these speculations let us turn to the proofs of the lemmas, beginning with
Lemma C. Throughout we will use the notation set up before the combinatorial version
of Theorem 2.1. Additionally, let us make two notes about terminology used throughout
the three lemmas. First, by ``natural" maps, we will mean certain maps derived from
the Weyl filtrations analogous to Pieri rules in [AB1]. These are the maps with which
we have to compose to get the appropriate maps between Hom groups as indicated in the
statements of the lemmas. Secondly, note the following abuse of notation. $K_\nu \otimes F$
is isomorphic to the skew Weyl module corresponding to the skew partition
$$(a_1+1)^{p_1} (a_2+1)^{p_2}\ldots (a_k+1)^{p_k} 1 / 1^{p_1 + \ldots + p_k} ,$$
so for simplicity we will call this skew partition $\nu \otimes 1$. For example we will
speak of tablueax of shape $\nu \otimes 1$, the generalized symmetrizer map
$d'_{\nu \otimes 1}$, etc. In the proof of Lemma C, often we will even write
$K_{\nu \otimes 1}$ meaning, of course, $K_\nu \otimes F$.

\vskip .3cm

{\it Proof of Lemma C.}
We will explicitly find a generator $f$ of Hom$(K_\la, K_{\nu \otimes 1})$, follow the
action of $f$ with the natural surjection $K_{\nu \otimes 1} \arr \hskip -.2cm \arr K_\la$,
and see which multiple of the identity we get in Hom$(K_\la, K_\la)$.

\vskip .3cm

Recall the results regarding skew Weyl modules summarized in [ABW, Theorem II.3.16].
Following the notation there we have $d'_\la: D_\la \arr \hskip -.2cm \arr K_\la$
with relations in $D_\la$ corresponding to each pair of adjacent rows in $\la$. So
$f$ (i.e., the map we seek) comes from a map $g: D_\la \arr  K_{\nu \otimes 1}$ that
sends these relations to 0. We will write a formula for a general map
$g: D_\la \arr  K_{\nu \otimes 1}$ and solve the resulting constraints on $g$ to
find $f$. In these calculations, we will have to rely very heavily on the material
sketched on pp.234-6 in [ABW]. More specifically, we will need (1) the definition of
the ``box" map to get relations defining $K_\la$ (these are described in general in
Note 2 at the end of Example 2 below), and (2) the procedure to ``straighten tableaux"
in $K_{\nu \otimes 1}$, which in turn is based on the box map associated to
$K_{\nu \otimes 1}$. (The second item is necessary as a systematic way to {\it find}
the answers but not necessary to verify them. However, straightening will play a more
essential role in the proof of Lemma A.) To undertand this material one should
additionally look at the careful treatment of the same topics in the dual case of
``Schur modules" on pp.226-232 of the same paper (statements of II.2.7 through II.2.16
and the proof of II.2.15). {\bf While following all subsequent calculations involving
tableaux it will be helpful to write the algebraic expressions pictorially in tableau
form.}

\vskip .3cm

Label an ordered basis of $F$ as $e_1,e_2,\ldots$ . Recall that for us standard
tableaux (called co-standard in [ABW]) will be fillings of Young diagrams with entries
from the chosen basis of $F$ such that the entries increase weakly along rows and
strictly along columns. {\it In what follows, we will for convenience identify any standard
tableau with the corresponding element of the appropriate tensor product of divided
powers of $F$.} For example, the canonical tableau (i.e., the tableau with all $e_1$'s
in the first row, all $e_2$'s in the second row and so on) $C_\la$ of shape $\la$ is
identified with $e_1^{(\la_1)} \otimes e_2^{(\la_2)} \otimes \ldots \in D_\la$. By the
method of weight submodules [AB2, pp.177-8], any map from $D_\la$ is specified completely
by its action on $C_\la$ (since $C_\la$ generates $D_\la$ as a ${\rm GL}(F)$-module) and such
a map may send $C_\la$ arbitrarily into the $\la$-weight submodule of the target module.
Together with the ``standard basis" for $K_{\nu \otimes 1}$ obtained from [ABW, Theorem
II.3.16], this gives us the following basis for Hom$(D_\la, K_{\nu \otimes 1})$.
$$\Bigl\{C_\la \mapsto d'_{\nu \otimes 1}(T_{\underline i})\; \big | \; T_{\underline i}
\hbox{ a standard tableau of shape } \nu \otimes 1 \hbox{ and weight } \la \Bigr\} .$$
Thus a general map $g: D_\la \arr  K_{\nu \otimes 1}$ may be written as
$C_\la \mapsto \sum_{\underline i} c_{\underline i} d'_{\nu \otimes 1}(T_{\underline i})$
with arbitrary integers $c_{\underline i}$. It is not very hard to explicitly describe
all the tableaux $T_{\underline i}$. After that one can explicitly write formulas for
the maps $\{C_\la \mapsto d'_{\nu \otimes 1}(T_{\underline i})\}$ in the above basis
by using appropriate polarizations (i.e., comultiplication followed by multiplication)
in the divided power algebra of $F$ (see examples below). It is now straightforward
but tedious to get $f$ following the strategy explained near the beginning of the
previous paragraph. Let us first work out three simple examples that will
simultaneously illustrate all the main features of the cases that we will have
to consider in the general proof.

\vskip .3cm

E{\sevenrm XAMPLE} 1. $\nu = a$, a partition with a single row. So $\la = a1$ and
$K_{\nu \otimes 1} = D_a \otimes F$. We have to find the smallest multiple
of identity in Hom$(K_{a1}, K_{a1})$ that factors through $D_a \otimes F$.
Note that here $d_{a1}': D_{a1} \arr K_{a1}$ with a single relation in
$D_{a1} = D_a \otimes F$ defining $K_{a1}$ (see below). To find a
generator $f$ of Hom$(K_{a1}, D_a \otimes F)$, let us look at the general
map $g$ in Hom$(D_a \otimes F, D_a \otimes F) \simeq a1$-weight submodule
of $D_a \otimes F$. The following two standard tableaux of shape $a \otimes 1$
form a basis for this weight submodule: $e_1^{(a)} \otimes e_2$ (i.e., the
canonical tableau $C_{a \otimes 1}$) and $e_1^{(a-1)}e_2 \otimes e_1$. These
correspond respectively to the maps $g_0$ = identity and $g_1$ = the
polarization defined by the composite map
$$D_a \otimes F \buildrel \Delta \over\larr D_{a-1} \otimes F \otimes F
\buildrel {m_{13}} \over \larr D_a \otimes F ,$$
where $\Delta$ = appropriate component of comultiplication on the first tensor factor
and $m_{13}$ = multiplication of the first and third tensor factors tensored with
identity on the second tensor factor. Now the map
$g = c_0 g_0 + c_1 g_1: D_a \otimes F \arr D_a \otimes F$
descends to a map from $K_{a1}$ exactly when it sends the kernel of
$d_{a1}': D_a \otimes F \arr K_{a1}$ to zero. This kernel is the image of
the comultiplication map $\Delta: D_{a+1} \arr D_a \otimes F$ and so is
generated by $\Delta$(the canonical tableau in $D_{a+1}$) =
$\Delta(e_1^{(a+1)}) = e_1^{(a)} \otimes e_1$. Evaluating $g$ on this, we get
$g(e_1^{(a)} \otimes e_1) = (c_0 + a c_1)(e_1^{(a)} \otimes e_1)$. Note the
coefficient $a$. It arises while calculating $m_{13}$ when we multiply
$e_1^{(a-1)}$ and $e_1$ in the divided power algebra of $F$. We now conclude
that $f$ is obtained by taking $c_0 = a$ and $c_1 = -1$. It remains to
calculate $d_{a1}' \circ f$ and see which multiple of the identity we get
in Hom$(K_{a1}, K_{a1})$. This is easiest to do by tracing what happens
to the canonical tableau. Applying $d_{a1}'$ to the images of the canonical
tableau $C_{a1}$ under $g_0$ and $g_1$ (i.e., to $e_1^{(a)} \otimes e_2$ and
$e_1^{(a-1)}e_2 \otimes e_1$), the desired integer is easily checked to be $a+1$.
The only thing to note is that
$d_{a1}'(e_1^{(a-1)}e_2 \otimes e_1) = - d_{a1}'(e_1^{(a)} \otimes e_2)$
by direct calculation.

\vskip .3cm

E{\sevenrm XAMPLE} 2. $\nu = a^2$, a partition with two rows of equal length $a$.
So $\la = a^21$ and $K_{\nu \otimes 1} = K_{a^2} \otimes F$. We have to
find the smallest multiple of identity in Hom$(K_{a^21}, K_{a^21})$ that
factors through $K_{a^2} \otimes F$.
To find a generator $f$ of Hom$(K_{a^21}, K_{a^2} \otimes F)$, let us look
at the general map $g$ in Hom$(D_{a^21}, K_{a^2} \otimes F) \simeq a^21$-weight
submodule of $K_{a^2} \otimes F$. A basis for this weight submodule is given by
images under $d_{a^2 \otimes 1}'$ of the three standard tableaux of shape $a^2 \otimes 1$
and weight $a^21$. Let us explicitly write this basis and the corresponding maps.

$d_{a^2 \otimes 1}'(e_1^{(a)} \otimes e_2^{(a)} \otimes e_3)$ corresponds to
$g_0 = d_{a^2 \otimes 1}'$.

$d_{a^2 \otimes 1}'(e_1^{(a-1)}e_2 \otimes e_2^{(a-1)}e_3 \otimes e_1)$ corresponds to
$g_1$ = the composite map
$$D_a \otimes D_a \otimes F \buildrel \Delta \over  \larr
D_{a-1} \otimes F \otimes D_{a-1} \otimes F \otimes F \larr
D_a \otimes D_a \otimes F \buildrel d_{a^2 \otimes 1}' \over \larr K_{a^2 \otimes 1} ,$$
where $\Delta$ is the appropriate comultiplication on the first
and second factors and the second map is $m_{14} \otimes m_{35} \otimes {\rm id}$,
i.e., multiplication on the indicated factors tensored with identity
on the second factor.

$d_{a^2 \otimes 1}'(e_1^{(a)} \otimes e_2^{(a-1)}e_3 \otimes e_2)$ corresponds to
$g_2$ = the composite map
$$D_a \otimes D_a \otimes F \buildrel \Delta \over \larr
D_a \otimes D_{a-1} \otimes F \otimes F \larr
D_a \otimes D_a \otimes F \buildrel d_{a^2 \otimes 1}' \over \larr K_{a^2 \otimes 1},$$
where $\Delta$ is the appropriate comultiplication on the second factor
and the second map is ${\rm id} \otimes m_{24} \otimes {\rm id}$, i.e., multiplication
on the indicated factors tensored with identity on the first and third factors.

\vskip 0.3cm

Now the map
$g = c_0 g_0 + c_1 g_1 + c_2 g_2: D_{a^21} \arr K_{a^2} \otimes F$
descends to a map from $K_{a^21}$ exactly when it sends the kernel of
$d_{a^21}': D_{a^21} \arr K_{a^21}$ to zero. This kernel is generated
by $a$ relations derived from the first two rows and a single relation
derived from the last two rows of the partition $a^21$. (See Note 2 at
the end of this example.) The former are images of the polarizations
involving first two rows
$$D_{a+t} \otimes D_{a-t} \otimes F \larr D_a \otimes D_{a} \otimes F$$
for $1 \leq t \leq a$, and so are generated by images of the canonical
tableau $e_1^{(a+t)} \otimes e_2^{(a-t)} \otimes e_3$, i.e., by
$e_1^{(a)} \otimes  e_1^{(t)}e_2^{(a-t)} \otimes e_3$ with $1 \leq t \leq a$.
The latter relation is similarly obtained by a polarization involving the last
two rows and is generated by $e_1^{(a)} \otimes  e_2^{(a)} \otimes e_2$.
Let us first use the last relation. By a calculation very similar to
the one in Example 1 and noting by direct calculation that the map $g_1$
already kills $e_1^{(a)} \otimes  e_2^{(a)} \otimes e_2$, we have
$$g\Big(e_1^{(a)} \otimes  e_2^{(a)} \otimes e_2 \Big) =
(c_0 + a c_2) d_{a^2 \otimes 1}'\Big(e_1^{(a)} \otimes e_2^{(a)} \otimes e_2 \Big) .$$
Now from the remaining $a$ relations, look at the one with $t = 1$.
Note that $g_0$ already kills  $e_1^{(a)} \otimes  e_1e_2^{(a-1)} \otimes e_3$. So
$$\eqalign{
g \Big(e_1^{(a)} \otimes  e_1e_2^{(a-1)} \otimes e_3 \Big) \; = \;
& c_1 d_{a^2 \otimes 1}'\Big(e_1^{(a-1)}e_2 \otimes e_1 e_2^{(a-2)} e_3 \otimes e_1 +
a \; e_1^{(a)} \otimes e_2^{(a-1)} e_3 \otimes e_1 \Big) + \cr
& c_2 d_{a^2 \otimes 1}'\Big(e_1^{(a)} \otimes e_2^{(a-1)} e_3 \otimes e_1 +
e_1^{(a)} \otimes e_1 e_2^{(a-2)} e_3 \otimes e_2 \Big) . \cr}
$$
Looking at the four terms inside two sets of parantheses on the right hand side,
the last one is killed by $d_{a^2 \otimes 1}'$ and the tableaux in middle two
terms are standard (in fact the same). For the first term, either by inspection
using the definition of $d_{a^2 \otimes 1}'$ or by using the ``straightening"
procedure in [ABW], we have
$$d_{a^2 \otimes 1}'\Big(e_1^{(a-1)}e_2 \otimes e_1 e_2^{(a-2)} e_3 \otimes e_1 \Big) \; =
\; - (a-1) d_{a^2 \otimes 1}'\Big(e_1^{(a)} \otimes e_2^{(a-1)} e_3 \otimes e_1 \Big) .$$
So we get
$$g \Big(e_1^{(a)} \otimes  e_1e_2^{(a-1)} \otimes e_3 \Big) \; = \;
(c_2 + c_1) d_{a^2 \otimes 1}' \Big(e_1^{(a)} \otimes e_2^{(a-1)} e_3 \otimes e_1 \Big) .$$
Equating the two relations treated above to zero, we get a one parameter family of
solutions for the integer coefficients $c_i$, namely $c_2 = -1, c_1 = 1$ and $c_0 = a$.
Since we are assured of the existence of $f$, and since ${\rm GL}(F)$-equivariance
of an integral multiple of any map clearly guarantees ${\rm GL}(F)$-equivariance
of the original map, the solution we found must give us the desired map $f$
generating Hom$(K_{a^21},K_{a^2}\otimes F)$. (So the remaining relations
must be automatically satisfied and were not necessary for our purpose.
This feature will repeat in the general calculation and can be understood
more conceptually as explained in Note 2 below.)

\vskip 0.3 cm

It remains to follow the action of $f$ with the natural surjection
$K_{a^2} \otimes F \arr K_{a^2 1}$
and see which multiple of the identity we get in Hom$(K_{a^2 1}, K_{a^2 1})$.
Applying $d_{a^2 1}'$ to the images of the canonical tableau $C_{a^2 1}$ under
$g_0, g_1$ and $g_2$ the desired integer is easily checked by direct calculation
to be $c_0 + c_1 - c_2 = a+2$.

\vskip 0.3 cm

{\it Notes.} (1) If $a=1$, the above calculation has to be modified since the terms where
$a-2$ occurs no longer make sense and have to be replaced by $0$. But it is easily
checked that the final answer $a+2=3$ is still valid.

\vskip 0.3 cm

(2) In general, to define any skew Weyl module by generators and relations, we need
for every pair of consecutive rows of the corresponding skew partition as many
relations as the length of the overlap between the two rows. Calling the lengths of
two such rows $p$ and $q$ and the length of their overlap $r$, the relations
corresponding to these rows are the images of the following $r$ polarizations
involving these two rows
$$\ldots D_{p+t} \otimes D_{q-t} \ldots \larr \ldots D_{p} \otimes D_{q} \ldots ,$$
where $q-r+1 \leq t \leq q$, ``$\ldots$" indicates tensoring by the divided powers
corresponding to the remaining rows and the above maps are identity on ``$\ldots$".
(Tracing the canonical tableau, it is easily seen from the definition of skew Weyl
modules that these are indeed relations, i.e., that $d'$ annihilates the images of
these maps. [ABW, Theorem II.3.16] shows, among other things, that these relations
suffice to define the skew Weyl module in question.) In fact it turns out that the
single relation with $t=q-r+1$ generates all the others up to a multiple (an argument
for this is sketched below) and so any equivariant map from the appropriate tensor
product of divided powers of $F$ satisfying this one relation must automatically satisfy
all the others.

\vskip 0.3 cm

Sketch for the claim in the previous sentence: letting the two rows be those numbered
$i-1$ and $i$, consider the map in ${\rm GL}(F)$ that takes $e_i$ to $e_{i-1} + e_i$ and
acts as the identity on other $e_j$. Apply the map induced on the divided power algebra
of $F$ to
$$ \ldots e_{i-1}^{(p)} \otimes e_{i-1}^{(q-r+1)}(e_{i-1} + e_i)^{(r-1)} \ldots ,$$
i.e., to the image of the canonical tableau under the polarization we picked out.
Expand the last divided power and notice that all the summands have different weights
and individually give all $r$ relations up to a multiple. (Also see lines 12-13 on
[ABW, p.209] but note the typo on line 13: $\la_2 - \mu_2 + 1$ should instead be
$\mu_1 - \mu_2 +1$.)

\vskip 0.3 cm

E{\sevenrm XAMPLE} 3. $\nu = ab$, a partition with two rows of unequal lengths $a>b$.
So $\la = ab1$ and $K_{\nu \otimes 1} = K_{ab} \otimes F$. We have to
find the smallest multiple of identity in Hom$(K_{ab1}, K_{ab1})$ that
factors through $K_{ab} \otimes F$. The outline of this calculation is
very similar to that in Example 2 with the following crucial change.
The relevant weight submodule now has dimension four instead of three,
resulting in significantly different constraints on the constants.
To find a generator $f$ of Hom$(K_{ab1}, K_{ab} \otimes F)$, let us look
at the general map $g$ in Hom$(D_{ab1}, K_{ab} \otimes F) \simeq ab1$-weight
submodule of $K_{ab} \otimes F$. A basis for this weight submodule is given by
images under $d_{ab \otimes 1}'$ of the four standard tableaux of shape $ab \otimes 1$
and weight $ab1$. Let us explicitly write this basis and the corresponding maps.
(The labeling of the maps in all three examples is chosen so as to be
consistent with the general case treated below.)

$d_{ab \otimes 1}'(e_1^{(a)} \otimes e_2^{(b)} \otimes e_3)$ corresponds to
$g_{00} = d_{ab \otimes 1}'$.

$d_{ab \otimes 1}'(e_1^{(a)} \otimes e_2^{(b-1)}e_3 \otimes e_2)$ corresponds to
$g_{01}$ = the composite map
$$D_a \otimes D_b \otimes F \buildrel \Delta \over \larr
D_a \otimes D_{b-1} \otimes F \otimes F \larr
D_a \otimes D_b \otimes F \buildrel d_{ab \otimes 1}' \over \larr K_{ab \otimes 1}  ,$$
where $\Delta$ is the appropriate comultiplication on the second factor and the
second map is ${\rm id} \otimes m_{24} \otimes {\rm id}$, i.e., multiplication on
the indicated factors tensored with identity on the first and third factors.

$d_{ab \otimes 1}'(e_1^{(a-1)}e_3 \otimes e_2^{(b)} \otimes e_1)$ corresponds to
$g_{10}$ = the composite map
$$D_a \otimes D_b \otimes F \buildrel \Delta \over \larr
D_{a-1} \otimes F \otimes D_b \otimes F \larr
D_a \otimes D_b \otimes F \buildrel d_{ab \otimes 1}' \over \larr K_{ab \otimes 1}  ,$$
where $\Delta$ is the appropriate comultiplication on the first factor and the
second map is $m_{14} \otimes s_{23}$, i.e., multiplication on the indicated
factors tensored with switching the middle two factors. (This is the new case
compared to Example 2, made possible by the fact that the first row
is longer than the second.)

$d_{ab \otimes 1}'(e_1^{(a-1)}e_2 \otimes e_2^{(b-1)}e_3 \otimes e_1)$ corresponds to
$g_{11}$ = the composite map
$$D_a \otimes D_b \otimes F \buildrel \Delta \over \larr
D_{a-1} \otimes F \otimes D_{b-1} \otimes F \otimes F \larr
D_a \otimes D_b \otimes F \buildrel d_{ab \otimes 1}' \over \larr K_{ab \otimes 1}  ,$$
where $\Delta$ is the appropriate comultiplication on the first
and second factors and the second map is $m_{14} \otimes m_{35} \otimes {\rm id}$,
i.e., multiplication on the indicated factors tensored with identity
on the second factor.

\vskip 0.3cm

Now the map
$g = c_{00} g_{00} + c_{01} g_{01} + c_{10} g_{10} + c_{11} g_{11}:
D_{ab1} \arr K_{ab} \otimes F$
descends to a map from $K_{ab1}$ exactly when it sends the kernel of
$d_{ab1}': D_{ab1} \arr K_{ab1}$ to zero. This kernel is generated
by $b$ relations derived from the first two rows and a single relation
derived from the last two rows of the partition $ab1$. The former
are generated by $e_1^{(a)} \otimes  e_1^{(t)}e_2^{(b-t)} \otimes e_3$ with
$1 \leq t \leq b$ and the latter by $e_1^{(a)} \otimes  e_2^{(b)} \otimes e_2$.
Let us first use the last relation. By a calculation similar to
the one in Examples 1 and 2, we have
$$g\Big(e_1^{(a)} \otimes  e_2^{(b)} \otimes e_2\Big) =
(c_{00} + b c_{01}) d_{ab \otimes 1}'\Big(e_1^{(a)} \otimes e_2^{(b)} \otimes e_2 \Big) +
(c_{10} + b c_{11}) d_{ab \otimes 1}'\Big(e_1^{(a-1)}e_2 \otimes e_2^{(b)} \otimes e_1 \Big)  .$$
Now from the remaining $a$ relations, look at the one with $t = 1$.
Note that $g_{00}$ already kills  $e_1^{(a)} \otimes  e_1e_2^{(b-1)} \otimes e_3$. So
$$\eqalign{
g \Big(e_1^{(a)} \otimes  e_1e_2^{(b-1)} \otimes e_3 \Big) \; = \;
& c_{01} d_{ab \otimes 1}'\Big(e_1^{(a)} \otimes e_2^{(b-1)} e_3 \otimes e_1 +
e_1^{(a)} \otimes e_1 e_2^{(b-2)} e_3 \otimes e_2 \Big) + \cr
& c_{10} d_{ab \otimes 1}'\Big(e_1^{(a-1)}e_3 \otimes e_1 e_2^{(b-1)} \otimes e_1 \Big) + \cr
& c_{11} d_{ab \otimes 1}'\Big(e_1^{(a-1)}e_2 \otimes e_1 e_2^{(b-2)} e_3 \otimes e_1 +
a \; e_1^{(a)} \otimes e_2^{(b-1)} e_3 \otimes e_1 \Big) . \cr
}$$
Looking at the five terms inside three sets of parantheses on the right hand side,
the second one is killed by $d_{ab \otimes 1}'$, the tableaux in the first and the
last terms are standard (in fact the same), and those in the third and fourth need to
be ``straightened." Either by inspection using the definition of $d_{ab \otimes 1}'$
or by using the straightening procedure in [ABW], for the third term we get
$$d_{ab \otimes 1}'\Big(e_1^{(a-1)}e_3 \otimes e_1 e_2^{(b-1)} \otimes e_1 \Big) \; = \;
- d_{ab \otimes 1}'\Big(e_1^{(a)} \otimes e_2^{(b-1)} e_3 \otimes e_1 \Big)  .$$
As for the fourth term, exactly as in Example 2, we have
$$d_{ab \otimes 1}'\Big(e_1^{(a-1)}e_2 \otimes e_1 e_2^{(b-2)} e_3 \otimes e_1 \Big) \; = \;
- (b-1) d_{ab \otimes 1}'\Big(e_1^{(a)} \otimes e_2^{(b-1)} e_3 \otimes e_1 \Big)  .$$
So we get
$$g\Big(e_1^{(a)} \otimes  e_1e_2^{(b-1)} \otimes e_3 \Big) \; = \;
\left( c_{01} - c_{10} + c_{11} (a - b + 1) \right)
d_{ab \otimes 1}'\Big(e_1^{(a)} \otimes e_2^{(b-1)} e_3 \otimes e_1 \Big)  .$$
Equating the two relations treated above to zero, we get a one parameter family of
solutions for the integer coefficients, namely $c_{11} = 1, c_{01} = -(a+1), c_{10} = -b$
and $c_{00} = (a+1)b$. As explained in Example 2, this must give us the desired
map $f$.

\vskip 0.3 cm

It remains to follow the action of $f$ with the natural surjection
$K_{ab} \otimes F \arr K_{ab1}$
and see which multiple of the identity we get in Hom$(K_{ab1}, K_{ab1})$.
Applying $d_{ab1}'$ to the images of the canonical tableau $C_{ab1}$ under
$g_{00}, g_{01}, g_{10}$ and $g_{11}$ the desired integer is easily checked
by direct calculation to be
$c_{00} - c_{01} - c_{10} + c_{11} = (a+1)b + (a+1) + b + 1 = (a+2)(b+1)$.
(Note that just as in Example 2, if $b = 1$ then terms involving $b-2$ have to be
dropped. But again the final result is easily checked to remain valid.)

\vskip 0.3 cm

Returning to the general case, the diagram of $\nu \otimes 1$ consists of $k$
rectangular blocks of rows (with the $j$-th block being $a_j^{p_j}$) plus a lone
box in the last row. Let us first describe the set of standard tableaux of shape
$\nu \otimes 1$ and weight $\la$. Such tableaux are in bijective correspondence
with $k$-tuples of integers ${\underline i}= i_1 \ldots i_k$ with $0 \leq i_j \leq p_j$.
For example, for $\nu = a_1^5 a_2^2 a_3^3$ , the index ${\underline i}= 302$
corresponds to the tableau with shape $\nu \otimes 1$, whose rows have rightmost
entries $e_1$ $e_2$ $e_4$ $e_5$ $e_9$ $e_6$ $e_7$ $e_8$ $e_{10}$ $e_{11}$ $e_3$
from top to bottom and whose all other entries match with the corresponding entries
in the canonical tableau $C_{\nu \otimes 1}$. The example will be clarified by the
general description given next.

\vskip .3cm

Let us describe $T_{\underline i}$ in general. It has the same entries as the
canonical tableau $C_{\nu \otimes 1}$ except possibly in the rightmost border
strip, in which the entries of $C_{\nu \otimes 1}$ undergo a cyclic permutation
as follows. When $i_j = 0$, the $j$-th block is unaffected, e.g., the second block
consisting of rows 6 and 7 in the above example. When the first nonzero $i_j$ from
the left appears, say $i_{j_1}$, the $i_{j_1}$-th entry in the rightmost column of
the $j_1$-th block is removed and the entries below it in the same block are moved
upward by one slot each. Now we find the next nonzero entry in ${\underline i}$, say
$i_{j_2}$, remove the $i_{j_2}$-th entry in the rightmost column of the $j_2$-th block
and place it in the empty space created at the bottom of the $j_1$ th block. Again the
entries below the removed one in the $j_2$-th block move one slot upward and so on,
until the last nonzero $i_{j_t}$ is used up. Now we remove the number in the last row
of $\nu \otimes 1$ and place it in the empty slot at the bottom of the $j_t$-th block.
Finally we place the very first entry we removed, i.e., the $i_{j_1}$-th entry in
the rightmost column of the $j_1$-th block, in the lone box in the last row of
$\nu \otimes 1$. Note for instance that $T_{\underline 0}$ = the canonical tableau
$C_{\nu \otimes 1}$. (It is a combinatorial exercise to prove that $T_{\underline i}$
are indeed all the standard tableaux of shape $\nu \otimes 1$ and weight $\la$.
But it is not really necessary to check this as long as the listed tableaux
suffice to produce a nontrivial map via our procedure.)

\vskip .3cm

The map $g_{\underline i}: D_\la \arr K_{\nu \otimes 1}$ corresponding to
$T_{\underline i}$ is obtained in a completely analogous manner to the maps
seen in the examples above. More precisely, in the first step we split off via
comultiplication a degree one piece from every row in which $T_{\underline i}$
differs from the canonical tableau $C_{\nu \otimes 1}$. Then we multiply what is
left for each such row (numbered, say, $r$) by the degree one piece split off from
$s$-th row where $e_s$ is the the last entry of $r$-th row in $T_{\underline i}$.
Finally we apply $d_{\nu \otimes 1}'$. Note for instance that $g_{\underline 0}$
is just $d_{\nu \otimes 1}'$.

\vskip .3cm

Let us now turn to the relations in $D_\la$ defining $K_\la$. These will help us
solve for the unknowns $c_{i_1 \ldots i_k}$. In what follows, sometimes we will
have to temporarily fix values of some of the components in ${i_1 \ldots i_k}$.
For ease of notation, after explaining such a choice, we will often denote such
fixed components simply by writing ``$\ldots$". The intended meaning will be clear
from the context. The relations defining $K_\la$ can be divided into three types
as follows.

\vskip .3cm

{\it Type} 1. The relation involving the last two rows of $\la$, i.e. the last row of the
partition $\nu$ and the last row of $\la$ consisting of exactly one box. Suppose
these rows are numbered $x$ and $x+1$ (so $x =  p_1 + \ldots + p_k$). Then the
relation in question is generated by $R_x = \ldots \otimes e_x^{(a_k)} \otimes e_x$,
where $\ldots$ indicates $x-1$ factors matching the first $x-1$ factors of the
canonical tableau $C_\la$. Applying $g$ to this relation, it is clear by explicit
calculation that unless $i_k = 0$ or $p_k$, $g_{\underline i}(R_x) = 0$, essentially
because otherwise all $a_k +1$ occurences of $e_x$ in $R_x$ are sent inside the last
block (which has only $a_k$ columns) by $g_{\underline i}$. (In Example 2 this was
manifested in the fact that the $g_1$ there already killed the relation in question.)
So let $i_k = 0$ or $p_k$ henceforth in this paragraph. For the moment arbitrarily
fix values of all other $i$'s so that we are considering only two of the maps
$g_{\underline i}$. Applying the corresponding two terms in $g$ to $R_x$, we get
exactly as in Examples 1 and 2 the expression
$(c_{\ldots 0} + a_k c_{\ldots p_k}) d_{\nu \otimes 1}'(S)$,
where $S$ is the standard tableau of shape $\nu \otimes 1$ obtained as follows.
Take $T_{\underline i}$ corresponding to either of the two $g_{\underline i}$ being
evaluated and obtain $S$ from $T_{\underline i}$ by replacing the single occurence
of the entry $e_{x+1}$ with $e_x$. Now allowing all possible choices of
${i_1 \ldots i_{k-1}}$, it is easy to see that the standard tableaux $S$ that
arise in the way just explained are all distinct. (Compare Example 3.) Altogether,
the relation under consideration gives us the following constraints.
$$c_{i_1 \ldots i_{k-1} 0} + a_k c_{i_1 \ldots i_{k-1} p_k} = 0 .$$
\vskip 0.3cm

{\it Type} 2. Relations involving consecutive rows in the same block. Let us see that this
case leads to a calculation virtually identical to the one in Example 2 where
a relation involving the first two rows was treated. Suppose we are dealing with
relations involving rows numbered $x$ and $x+1$ and that these are respectively
the $y$-th and $(y+1)$-th rows of the $j$-th block (and so each is of length $a_j$).
Then one of the relevant relations is generated by
$R_x = \ldots \otimes  e_x^{(a_j)} \otimes e_x e_{x+1}^{(a_j-1)} \otimes \ldots$ ,
where $\ldots$ indicates expressions identical to the corresponding factors of the
canonical tableau $C_\la$. Applying $g$ to this relation results in the following.
For $i_j$ other than $y$ and $y+1$, $g_{\underline i}(R_x) = 0$ by explicit
calculation, essentially because all $a_j+1$ occurences of $e_x$ in $R_x$ are kept
within the $j$-th block of $\la$ (which has only $a_j$ columns) by such
$g_{\underline i}$. (In Example 2 this was manifested in the fact that the $g_0$
there already killed the relation in question.) So let $i_j = y$ or $y+1$ henceforth
in this paragraph. For the moment arbitrarily fix values of all other $i$'s so that
we are considering only two of the $g_{\underline i}$. Then exactly as in Example 2,
we get an expression with four terms inside two sets of parantheses with obvious
changes in subscripts and placement of $e$'s. By a very similar calculation, this
expression simplifies to
$(c_{\ldots y \ldots} + c_{\ldots (y+1) \ldots})d_{\nu \otimes 1}'(T)$, where $y$
and $y+1$ appear in the $j$-th slot and $T$ is the following standard tableau of
shape $\nu \otimes 1$. From the two $g_{\underline i}$ being evaluated, consider the
one with $i_j = y$ and take the corresponding tableau $T_{\underline i}$. Obtain $T$
from this $T_{\underline i}$ by replacing the single $e_{x+1}$ in the $x$-th row
with $e_x$. It is now easy to see that the tableaux $T$ that arise in this fashion
are all distinct for distinct choices of ${\underline i}$. (Of course only $k-1$
components of ${\underline i}$ are being chosen here, in view of our earlier reasoning
for $i_j$.) So just as we got $c_1 + c_2 = 0$ in Example 2, we get the constraint
$c_{\ldots y \ldots} + c_{\ldots (y+1) \ldots } = 0$, where $y$ and $y+1$ appear in
the $j$-th slot. It should be clear now that in general such relations give us the
following equations.
$$c_{i_1 \ldots i_j \ldots i_k} + c_{i_1 \ldots (i_j+1) \ldots i_k} = 0 ,
\quad 1 \leq j \leq k , \quad 1 \leq i_j \leq p_j -1 .$$
In view of these relations it is enough to determine $c_{\underline i}$ with
$i_j = 0$ or $p_j$ for each $j$.

\vskip .3cm

{\it Type} 3. Relations involving the last row of one block and the first in the next. This
case leads to calculation very similar to the corresponding calculation in Example
3. Suppose we are dealing with relations involving rows numbered $x$ and $x+1$ and
that these are respectively the last row of the $j$-th block and the first row of
the $(j+1)$-th block. So row $x$ is of length $a_j$ and row $x+1$ of smaller length
$a_{j+1}$. Then one of the relevant relations is generated by
$R_x = \ldots \otimes e_x^{(a_j)} \otimes e_x e_{x+1}^{(a_{j+1}-1)} \otimes \ldots$ ,
where $\ldots$ indicates expressions identical to the corresponding factors of the
canonical tableau $C_\la$. Consider the calculation of $g(R_x)$. Note that all
$a_j + 1$ occurences of $e_x$ in $R_x$ are contained within the $j$-th and $(j+1)$-th
blocks, where the number of columns is never more than $a_j$. So for any
$g_{\underline i}$ that leaves this property unchanged, clearly
$g_{\underline i}(R_x) = 0$. Using this it is easy to see that for $i_j$ other than
$0$ and $p_j$, $g_{\underline i}(R_x) = 0$ and furthermore, if $i_j = 0$ then again
$g_{\underline i}(R_x) = 0$ unless $i_{j+1} = 1$. Thus the terms involving only those
$g_{\underline i}$ where either ($i_j = 0$ and $i_{j+1} = 1$) or ($i_j = p_j$ and
$i_{j+1} = 0, \ldots, p_{j+1}$) can possibly survive. Until further notice arbitrarily
fix values of all $i$'s other than $i_j$ and $i_{j+1}$. So in view of the preceding
discussion, we are considering only $p_{j+1} + 2$ terms. Let us consider these one by
one after setting up some more notation. $T_{\underline i}^x$ will denote the (nonstandard)
tableau obtained from $T_{\underline i}$ by replacing the first appearance of $e_{x+1}$
in row $x+1$ with $e_x$. (See Remark 1 at the end of this proof  if there is no such
appearance.) Define a standard tableau $U$ as follows. Consider the $T_{\underline i}$
with $i_j = 0$ and $i_{j+1} = 1$. Obtain $U$ by replacing the only occurence of $e_{x+1}$
outside row $x+1$ with $e_x$.

\vskip .3cm

(i) For the ${\underline i}$ with $i_j = 0$ and $i_{j+1} = 1$, $g_{\underline i}(R_x)
= d_{\nu \otimes 1}'(T_{\underline i}^x + U) = d_{\nu \otimes 1}'(U)$ by noting that
$d_{\nu \otimes 1}'(T_{\underline i}^x) = 0$. (Compare the very similar evaluation of
$g_{01}$ in Example 3 on the relevant relation.)

\vskip .3cm

(ii) For the ${\underline i}$ with $i_j = p_j$ and $i_{j+1} = 0$, $g_{\underline i}(R_x)
= d_{\nu \otimes 1}'(T_{\underline i}^x) = (-1)^{p_{j+1}}d_{\nu \otimes 1}'(U)$. The
second of these equalities is obtained by easy application of the ``straightening"
procedure in [ABW]
or checked even more easily by explicit calculation using the definition of
$d_{\nu \otimes 1}'$. (Compare the very similar evaluation of $g_{10}$ in Example 3
on the relevant relation.)

\vskip .3cm

(iii) For the ${\underline i}$ with $i_j = p_j$ and $i_{j+1} = 1$,
$$g_{\underline i}(R_x) = d_{\nu \otimes 1}'(T_{\underline i}^x + a_j U) =
(a_j - a_{j+1} + 1) d_{\nu \otimes 1}'(U)$$
since by straightening or explicit calculation,
$d_{\nu \otimes 1}'(T_{\underline i}^x) = -(a_{j+1} - 1)d_{\nu \otimes 1}'(U)$.
(Compare the very similar evaluation of $g_{11}$ in Example 3 on the relevant
relation.)

\vskip .3cm

(iv) For the ${\underline i}$ with $i_j = p_j$ and $i_{j+1} = r, 2 \leq r \leq p_{j+1}$,
$g_{\underline i}(R_x) = d_{\nu \otimes 1}'(T_{\underline i}^x)
= (-1)^{r-1} d_{\nu \otimes 1}'(U)$ by straightening or explicit calculation. (This
case does not have a counterpart in Example 3, since there $p_{j+1}$ was 1. But the
necessary simplification follows the same pattern as in item (ii) above.)

\vskip .3cm

Collecting the information in items (i) through (iv) and using the constraints obtained
from relations of type 2 to convert all $i_{j+1}$ into $0$ or $p_{j+1}$, we get
$$g(R_x) = (-1)^{p_{j+1} -1} \Big( c_{\ldots 0 p_{j+1} \ldots} - c_{\ldots p_j 0 \ldots} +
(a_j - a_{j+1} + p_{j+1}) c_{\ldots p_j p_{j+1} \ldots} \Big) d_{\nu \otimes 1}'(U) ,$$
where only the $j$-th and $(j+1)$-th components of ${\underline i}$ are written.
Now allowing ${\underline i}$ to vary, it is easy to see that the standard tableaux
$U$ that arise as explained above are all distinct. So we get the following
constraints, where again only the $j$-th and $(j+1)$-th components in ${\underline i}$
are written, all others being assumed to be the same in each term of the equation.
For $ 1 \leq j \leq k-1$,
$$c_{0 p_{j+1}} - c_{p_j 0} + (a_j - a_{j+1} + p_{j+1}) c_{p_j p_{j+1}} = 0 .$$
Let us now solve all the constraints on $c_{\underline i}$ obtained above.
Order all $c_{\underline i}$ lexicographically with respect to the indices
${\underline i}$. (In each constraint above the variables are written in
lexicographic order.) It is easily seen that the constraints give us a
homogeneous linear system of integer equations whose coefficient matrix is
triangular with all 1's on the diagonal and that $c_{p_1 \ldots p_k}$ is the
only free variable. So we have a one-parameter family of solutions for
$c_{\underline i}$ and each solution is uniquely specified by assigning an
arbitrary integer value to $c_{p_1 \ldots p_k}$. We will verify that
$$b_{\underline i}
\; = \; \epsilon_{\underline i}  \prod_{i_j = 0} (a_j+ p_{j+1} + p_{j+2} + \ldots + p_k)
\; = \; \epsilon_{\underline i}  \prod_{i_j = 0} (h_j - p_j)
$$
is a generator of this family, where $\epsilon_{\underline i}$ is the sign of the
cyclic permutation in the rightmost border strip of $T_{\underline i}$ with reference
to the canonical tableau $C_{\nu \otimes 1}$. Since
$\epsilon_{i_1 \ldots i_{k-1} 0} = - \epsilon_{i_1 \ldots i_{k-1} p_k}$ and
$b_{i_1 \ldots i_{k-1} 0} = a_k b_{i_1 \ldots i_{k-1} p_k}$ the constraints of type
1 are satisfied. The sign $\epsilon_{\underline i}$ ensures that the constraints of
type 2 are satisfied. As for the constraints of type 3, noting that the sign of the
$b_{\underline i}$ corresponding to the middle term of this constraint is opposite
to that for the first and the third term, the verification boils down to the simple
fact that
$$(a_j+ p_{j+1} + p_{j+2} + \ldots + p_k) - (a_j  - a_{j+1} + p_{j+1}) -
(a_{j+1} + p_{j+2} + \ldots + p_k) = 0 .$$

Since $|b_{p_1 \ldots p_k}| = 1$, we have the desired map $f$. After having computed
$f$, we need to follow its action with the natural surjection $K_{\nu \otimes 1} \arr K_\la$.
It is easy to see that under this surjection $d_{\nu \otimes 1}'(T_{\underline i})$
is sent to $\epsilon_{\underline i} d_{\la}'(C_{\la})$. Therefore the map
Hom$(K_\la, K_{\nu \otimes 1}) \arr {\rm Hom} (K_\la, K_\la)$ is given by the integer
$\sum_{\underline i} |b_{\underline i}|$. Call this expression $E_\nu$. We will show by
induction on the number of blocks $k$ in the partition $\nu$ that $E_\nu = h_1 \ldots h_k$.
The base case $k=1$ is an easy check (or even subsumed in the induction by defining
$E_\phi = 1$ and seeing that this makes obvious sense in the context of the lemma).
Let $\nu' = {a_2^{p_2}\ldots a_k^{p_k}}$ be the partition obtained from $\nu$ by deleting
the first block $a_1^{p_1}$. The terms in $E_\nu$ with $i_1 = 0$ give $(h_1-p_1) E_{\nu'}$.
The remaining terms in $E_\nu$, by fixing $i_1 = 1, \ldots, p_1$ at a time, can be
partitioned into $p_1$ groups each of which is $E_{\nu'}$. This completes the proof
of Lemma C.

\vskip .3cm

{\it Remarks.} (1) If $a_k = 1$, one has to modify a few of the details in the calculations
required to find the constraints, but the same constraints stay valid. Specifically,
while dealing with relations of type 2 for rows in the last block, two of the four
terms within two sets of parantheses have to replaced by 0 (see Note 1 at the end of
Example 2). While dealing with relations of type 3, the definition of tableaux
$T_{\underline i}^x$ no longer makes sense when $j+1=k$ and $i_k = 1$ because then
the single entry in row $x+1$ of $T_{\underline i}$ is not $e_{x+1}$. But this
simply means that the $T_{\underline i}^x$ in items (i) and (iii) has to be replaced
by 0, leaving the end result in these items unchanged for the case under consideration
(see the Note at the end of Example 3).

\vskip .3cm

(2) The map $f$ may be of independent interest because in a sense it gives us an explicit
formula for the intertwining homomorphism between neighboring Weyl modules for ${\rm GL}_n$,
see [Andersen1, Section 6].

\vskip .3cm

In proofs of the remaining lemmas, we will need several items from the proof of Lemma C.
These include the ordered basis $e_1, e_2, \ldots$ for $F$, the indices ${\underline i}$,
the tableaux $T_{\underline i}$ and the associated maps $g_{\underline i}$, the numbers
$b_{\underline i}$ and the map $f$.

\vskip .3cm

{\it Proof of Lemma B.}
Let us first identify the module $N$ in the statement of the lemma as a skew Weyl
module. Using [AB1, 6.6] and contravariant duality, we have the short exact sequence
$0 \arr K_\mu \arr K_{\nu} \otimes F \buildrel \pi \over \arr K_\xi \arr 0$ ,
where $K_\xi$ is the skew Weyl module corresponding to the  skew partition
$\xi  = a_1^{p_1 + 1} a_2^{p_2} \ldots a_k^{p_k} / a_1 - 1$. (In words, the diagram
of $\xi$ is obtained by placing a single box immediately above the last box in the
first row of $\nu$.) So $N = K_\xi$ and our task is to find the smallest multiple of
identity in Hom$(K_\la, K_\la)$ that factors through $K_\xi$. For this we will find
a generator $f'$ of Hom$(K_\la, K_\xi)$ and follow its action with the natural
surjection $K_\xi \arr K_\la$ (obtained from the natural surjection
$K_\nu \otimes F \arr K_\la$ used in Lemma C, which kills the submodule $K_\mu$. See
the exact sequence above.) The generator $f'$ can be found exactly as $f$ was
found in Lemma C. But instead we will use the work already done to find $f$ and finish
the proof as follows. Consider the composite map
$\pi \circ f: K_\la \arr K_\xi$ where $\pi$ is the surjection in the exact sequence
above. Using the standard basis for $K_\xi$, we will write an explicit expression for
$\pi \circ f(d_\la'(C_\la))$. It will be obvious that this expression is divisible
by $h_1$ and by no larger integer. So $f'= \pi \circ f/h_1$ is a generator for
Hom$(K_\la, K_\xi)$ and composing it with the natural surjection $K_\xi \arr K_\la$
gives the desired result in view of Lemma C.

\vskip .3cm

Let us lay some groundwork for the calculation of $\pi \circ f(d_\la'(C_\la))$. First note
that removing the lone box in the last row of $\nu \otimes 1$ and placing it directly above
the last box in the first row of $\nu$ gives us the diagram of the skew partition $\xi$.
We will need this relocation on two occasions. In the first instance it is involved
in the surjection $\pi : K_\nu \otimes F \arr K_\xi$. Explicitly $\pi$ maps
$d_{\nu \otimes 1}'(C_{\nu \otimes 1})$ to
$d_{\nu'}'(C_ {\nu'})\otimes (e_m \wedge e_1 \wedge \ldots \wedge e_{p_1})$, where
$\nu'$ is the partition obtained by stripping off the last column of $\nu$ and $e_m$
is the entry in the last row of $C_ {\nu \otimes 1}$. (All this is only a notational
issue arising from the simple fact that the diagram of $\nu$ and a single box can be
``joined cornerwise" in either order to give us two skew partitions whose associated
skew Weyl modules are both isomorphic to $K_\nu \otimes F$.) Secondly the same
relocation is involved in describing the standard tableaux of shape $\xi$ and weight
$\la$. These are $T_{i_2 i_3 \ldots i_k}'$, defined as follows. $i_2 i_3 \ldots i_k$
is an arbitrary $(k-1)$-tuple of integers satisfying $0 \leq i_j \leq p_j$ and
$T_{i_2 i_3 \ldots i_k}'$ is identical to the tableau $T_{1i_2 i_3 \ldots i_k}$ in
the proof of Lemma C except for the relocation of a single box (with its content intact)
required to get the diagram of $\xi$ from that of $\nu \otimes 1$. (Again it is an
optional combinatorial exercise to verify that these are precisely the tableaux of the
given description. Alternatively this fact is immediate from the corresponding
combinatorial statement in the proof of Lemma C.) Using all the setup in this paragraph,
it is an easy check that $\pi$ maps $d_{\nu \otimes 1}'(T_{i_1 i_2 \ldots i_k})$
to $(-1)^{i_1 -1} d_{\xi}'(T_{i_2 \ldots i_k}')$ if $1 \leq i_1 \leq p_1$ and to
$(-1)^{p_1} d_{\xi}'(T_{i_2 \ldots i_k}')$ if $i_1 = 0$. Putting everything together
$$\pi \circ f(d_\la'(C_\la))
= \sum_{i_2 \ldots i_k} \Big( (-1)^{p_1} b_{0 i_2 \ldots i_k} + \sum_{i_1=1}^{p_1}
                                    (-1)^{i_1-1} b_{i_1 i_2 \ldots i_k}
                              \Big) d_{\xi}'(T_{i_2 \ldots i_k}') . $$
Referring to the values of $b_{\underline i}$, we have
$b_{0 i_2 \ldots i_k} = - (h_1 - p_1) b_{p_1 i_2 \ldots i_k}$ and
$(-1)^{i_1-1} b_{i_1 i_2 \ldots i_k} = (-1)^{p_1 -1} b_{p_1 i_2 \ldots i_k}$.
So the coefficient of $d_{\xi}'(T_{i_2 i_3 \ldots i_k}')$ in the above expression works
out to be $(-1)^{p_1 -1} h_1 b_{p_1 i_2 \ldots i_k}$. Since $|b_{p_1 p_2 \ldots p_k}| = 1$,
the proof of Lemma B is complete.

\vskip .3cm

{\it Proof of Lemma A.}
Our task is to find the smallest multiple of the identity in Hom$(K_\nu, K_\nu)$ that
factors through $K_{\la/1}$. For this we will compose the natural injection
$\iota: K_\nu \hookrightarrow K_{\la/1}$ with a generator of Hom$(K_{\la/1}, K_\nu)$.
This generator can be found exactly as in Lemma C. But instead we will again get it
using the work already done for Lemma C by employing an idea from [AB1] as follows.

\vskip .3cm

In this paragraph only, consider $K_\la, K_{\la/1}$ and $K_\nu$ as functors. Applied to
the free abelian group $F$ these give us the representations $K_\la(F), K_{\la/1}(F)$
and $K_\nu(F)$ of ${\rm GL}(F)$. Consider an extra copy of {\bf Z} (which we will regard as
the trivial representation of ${\rm GL}(F))$ with basis $e_0$ and use the ordered basis
$e_0, e_1, \ldots$ for ${\bf Z} \oplus F$. Now the module $K_{\la}({\bf Z} \oplus F)$
considered as a representation of ${\rm GL}(F)$ splits into a direct sum by the content of
$e_0$. Using the standard basis theorem, the summand with $e_0$-content one is spanned
by standard tableaux with a single $e_0$ entry in the top left corner of $\la$ and
is clearly isomorphic to $K_{\la/1}(F)$. Consider the following composite map of
representations of ${\rm GL}(F)$ beginning with the containment just explained.
$$   K_{\la/1}(F) \arr K_{\la}({\bf Z} \oplus F)
\buildrel f \over \arr K_{\nu}({\bf Z} \oplus F) \otimes ({\bf Z} \oplus F)
\buildrel pr \over \arr K_\nu({\bf Z} \oplus F) ,$$
where the second map is the map $f$ in the proof of Lemma C (now between representations
of ${\rm GL}({\bf Z} \oplus F)$) and the third map is the projection $pr$ derived from the
direct sum decomposition $K_\nu({\bf Z} \oplus F) \otimes ({\bf Z} \oplus F) \simeq
K_\nu({\bf Z} \oplus F) \oplus (K_\nu({\bf Z} \oplus F) \otimes F)$. Let us trace the
action of this composite map on the canonical tableau $C_{\la/1}$ of shape $\la/1$.
The first map takes $d_{\la/1}'(e_1^{(\la_1 -1)} \otimes e_2^{(\la_2)} \otimes \ldots)$
to $d_{\la}'(e_0 e_1^{(\la_1 -1)} \otimes e_2^{(\la_2)} \otimes \ldots)$. Now
using the work done for Lemma C,
$$f\left( d_{\la}'\Big(e_0 e_1^{(\la_1 -1)} \otimes e_2^{(\la_2)} \otimes \ldots\Big) \right) =
\sum_{\underline i} b_{\underline i} \, g_{\underline i}
\Big(e_0 e_1^{(\la_1 -1)} \otimes e_2^{(\la_2)} \otimes \ldots\Big) .$$
Applying the projection $pr$ to this element, we see that the terms involving
$g_{\underline i}$ with $i_1 \neq 1$ are killed. This is because applying such
$g_{\underline i}$ would result in a tableau having an entry other than $e_0$
in the lone box in the last row of $\nu \otimes 1$, so this tableau would be
killed by $pr$. As for
$g_{1 i_2 \ldots i_k} (e_0 e_1^{(\la_1 -1)} \otimes e_2^{(\la_2)} \otimes \ldots)$,
what survives upon projection is exactly $d'_{\nu}(T_{i_2 \ldots i_k}'')$, where
the standard tableau $T_{i_2 \ldots i_k}''$ of shape $\nu$ is the same as the tableau
$T_{1 i_2 \ldots i_k}$ in the proof of Lemma C (or $T_{i_2 \ldots i_k}'$ in the proof
of Lemma B) with the extra box containing the entry $e_1$ removed. (Again it is an
optional combinatorial check -- or immediate from the corresponding fact in proofs of
Lemma C or Lemma B -- that $T_{i_2 \ldots i_k}''$ exhaust the standard tableaux
of shape $\nu$ and weight $e_1^{\la_1 -1} e_2^{\la_2} e_3^{\la_3} \ldots$ .) It is also
clear from the calculation (or even a priori) that the image of the composite map above
lies in the submodule $K_{\nu}(F)$ of $K_{\nu}({\bf Z} \oplus F)$. All in all, we have
constructed a map $f'': K_{\la/1}(F) \arr K_\nu(F)$ such that
$$f''\Big( d_{\la/1}'\Big(e_1^{(\la_1 -1)} \otimes e_2^{(\la_2)} \otimes \ldots \Big) \Big) =
\sum_{\underline i = 1i_2 \ldots i_k} b_{1 i_2 \ldots i_k}  d'_{\nu}(T_{i_2 \ldots i_k}'') .$$
Recall how in the proof of Lemma C a formula for $f$ was written using maps
$g_{\underline i}$. In exactly the same fashion from the preceding equation we can say
that $f''$ descends from the map
$\sum_{1i_2 \ldots i_k} b_{1 i_2 \ldots i_k} g_{i_2 \ldots i_k}'': D_{\la/1}(F) \arr K_\nu(F),$
where $g_{i_2 \ldots i_k}''$ is the map ``corresponding" to the tableau $T_{i_2 \ldots i_k}''$.
(In other words $g_{i_2 \ldots i_k}''$ can be built from the tableau $T_{i_2 \ldots i_k}''$
using polarizations in the divided power algebra of $F$ exactly the way $g_{\underline i}$ was
built from the tableau $T_{\underline i}$ in the proof of Lemma C, see the paragraph immediately
after the one describing $T_{\underline i}$.) Further, since $|b_{1 p_1 \ldots p_k}| = 1$,
$f''$ is indivisible and hence a generator of Hom$(K_{\la/1}(F), K_\nu(F))$.

\vskip .3cm

Let us now compose the the map $f''$ in the previous paragraph with the natural
injection $\iota: K_\nu \hookrightarrow K_{\la/1}$. $\iota$ descends from a map that
``polarizes a degree one piece" from each row into the next row. Recall that this
means splitting off a degree one component from one row via diagonalization and
then multiplying this component into another row, both operations being done in the
divided power algebra of $F$. Let us make this explicit using the canonical
tableau $C_\nu = e_1^{(\nu_1)} \otimes e_2^{(\nu_2)} \otimes \ldots$ . We have
$\iota(d'_\nu(C_\nu)) = d'_{\la/1}
(e_1^{(\la_1 - 1)} \otimes e_1 e_2^{(\la_2 - 1)} \otimes e_2 e_3^{(\la_3 - 1)} \otimes \ldots)$.
(Note that $\la$ and $\nu$ are identical with the exception that $\la$ has an extra row
consisting of a single box.) We will find the image of this element under $f''$ by computing
the individual terms
$g_{i_2 \ldots i_k}''(e_1^{(\la_1 - 1)} \otimes e_1 e_2^{(\la_2 - 1)} \otimes
e_2 e_3^{(\la_3 - 1)} \otimes \ldots)$.
This computation involves two steps. (1)Polarize degree one pieces from several rows into
previous rows as dictated by the entries in the rightmost border strip of the tableau
$T_{i_2 \ldots i_k}''$. (2)Apply $d'_\nu$ to the result of the first step. Since the result
of the first step is a linear combination of tableaux that are in general not standard,
in the second step one has to straighten these tableaux using the procedure in [ABW]. To
keep control of the calculation it will be convenient for us to mix the order of operations
involved in steps 1 and 2 as follows. Proceeding from top row to the bottom row, we will
follow each {\it single} polarization immediately by straightening. Each such straightening
will involve only a fragment of a tableau up to the rows involved the preceding polarization.
It will be clear that the end result is unaffected by such interlacing of steps 1 and 2 used
for one polarization at a time.

\vskip 0.3cm

Let us illustrate the above discussion by first treating an extreme case. Consider
$g_{1 \ldots 1}''(e_1^{(\la_1 - 1)} \otimes e_1 e_2^{(\la_2 - 1)} \otimes
e_2 e_3^{(\la_3 - 1)} \otimes \ldots)$.
Note that $g_{1 \ldots 1}''$ involves, for each pair of consecutive nonzero rows of $\la$,
polarization of a degree one component from the lower row into the immediately preceding
row. For now look at what happens after doing only the polarization from the second row
into the first row, which affects only the first two tensor factors. This gives
$$e_1^{(\la_1 - 1)} \otimes e_1 e_2^{(\la_2 - 1)} \otimes  \ldots \; \longmapsto \;
\la_1  \Big(e_1^{(\la_1)} \otimes  e_2^{(\la_2 - 1)} \otimes  \ldots \Big) +
\Big(e_1^{(\la_1 - 1)} e_2 \otimes e_1 e_2^{(\la_2 - 2)} \otimes  \ldots \Big) .$$
Notice that the tableau fragment displayed in second term on the right hand side is
non-standard due to a violation in the very first column. It is clear by looking that
this non-standardness will persist after subsequent polarizations. So after the application
of $d'_\nu$, we will need to apply the straightening procedure to the first two rows.
(And lower rows too, but we will deal with that later. Here ``straightening" means
replacing a tableau by a linear combination of tableaux as prescribed in [ABW] that
will give the same result upon applying the appropriate generalized symmetrizer map $d'$.)
Moreover, we may do this straightening {\it before} applying the rest of the polarizations
involved in $g_{1 \ldots 1}''$ without affecting the overall result. This is because
subsequent polarizations will only involve rows numbered two and below. Rows below the
second are entirely unaffected by the results of the proposed straightening. As for the
second row, the next polarization will result in a degree one piece being multiplied into
it, but this multiplication and the proposed straightening together give the same result
for the second row regardless of the order in which they are performed. This is simply from
associativity of multiplication. To carry out the proposed straightening, consider the
diagonalization
$\Delta(e_1^{(\la_1)}  e_2) =  e_1^{(\la_1)} \otimes  e_2  +  e_1^{(\la_1 - 1)} e_2 \otimes e_1.$
Using this the nonstandard term straightens to
$ - (\la_2- 1) (e_1^{(\la_1)} \otimes  e_2^{(\la_2 - 1)} \otimes  \ldots)$, where the
constant $(\la_2 - 1)$ is due to the multiplication of $e_2$ and $e_2^{(\la_2 - 2)}$ in the
divided power algebra. Combining, the result so far in calculating $g_{1 \ldots 1}''$
can be shown as follows.
$$e_1^{(\la_1 - 1)} \otimes e_1 e_2^{(\la_2 - 1)} \otimes  \ldots \; \longmapsto \;
(\la_1  - \la_2  + 1) \Big(e_1^{(\la_1)} \otimes  e_2^{(\la_2 - 1)} \otimes  \ldots \Big) .$$
(Note that this is a just a schematic representation of what happens after
applying steps 1 and 2 for a single polarization. In particular we cannot write
$d'_\nu$ on the right hand side until all polarizations are applied.) Since the
fragment obtained so far matches that for the standard tableau $C_\nu$, clearly
subsequent polarizations in the calculation of $g_{1 \ldots 1}''$ will not result in
terms that are nonstandard in the first two rows. So evidently we may use the same
logic on successive pairs of rows of $\la$. Inductively we get the end result to be
$d'_\nu(C_\nu)$ times the product of $\la_{t-1} - \la_{t} + 1$ over successive pairs
of nonzero rows of $\la$.

\vskip 0.3cm

By an extension of the above argument, we will show that in general,
$$g_{i_2 \ldots i_k}''\Big(e_1^{(\la_1 - 1)} \otimes e_1 e_2^{(\la_2 - 1)} \otimes  \ldots\Big)
\; = \;
(-1)^{p_1 + \ldots + p_k} \, \epsilon_{1 i_2 \ldots i_k}
\prod_{t \in S} (\la_{t-1} - \la_{t} + 1) d'_\nu(C_\nu) , \leqno(*)$$
where the product is taken over the set $S$ of nonzero rows numbered $t$ such that
$g_{i_2 \ldots i_k}''$ involves a polarization of the $t$-th row of $\la$ into a previous
row of $\la$. (Recall from Lemma C that $\epsilon_{\underline i}$ is the sign of the cyclic
permutation in the rightmost border strip of $T_{\underline i}$ with reference to the canonical
tableau $C_{\nu \otimes 1}$ and that $p_1 + \ldots + p_k$ is the number of rows in $\nu$.)

\vskip 0.3cm

To prove the claim, fix an arbitrary $g_{i_2 \ldots i_k}''$ and as before, let us consider
just the first polarization involved in calculating it. Let us suppose that the first row
that it polarizes (necessarily into the first row) is that numbered $t$. Then we have
$$ \eqalign{
&e_1^{(\la_1 - 1)} \otimes e_1 e_2^{(\la_2 - 1)} \otimes  \ldots \otimes
e_{t-1} e_t^{(\la_t - 1)} \otimes \ldots  \quad \longmapsto \cr
&\hskip 3cm e_1^{(\la_1 - 1)} e_{t-1} \otimes e_1 e_2^{(\la_2 - 1)} \otimes  \ldots \otimes
e_t^{(\la_t - 1)} \otimes \ldots \; +  \cr
& \hskip 3cm e_1^{(\la_1 - 1)} e_{t} \otimes e_1 e_2^{(\la_2 - 1)} \otimes
\ldots \otimes e_{t-1} e_t^{(\la_t - 2)} \otimes \ldots . \cr
}$$
Now in general, after applying $d'_\nu$, both fragments displayed on the right hand side will
need to be straightened. To straighten the first term (needed if $t > 2$), we use
$\Delta(e_1^{(\la_1)} e_{t-1}) = e_1^{(\la_1)} \otimes e_{t-1} + e_1^{(\la_1 - 1)} e_{t-1} \otimes e_1.$
This results in
$- e_1^{(\la_1)} \otimes e_2^{(\la_2 - 1)}e_{t-1} \otimes  \ldots \otimes
e_t^{(\la_t - 1)} \otimes \ldots $,
i.e., the first term undergoes an exchange of $e_{t-1}$ and $e_1$ between the first and second
rows and picks up a negative sign. Now standardness is violated between second and third rows
(unless $t=3$). So repeat the same procedure using the second and third rows, and so on until
$e_{t-1}$ moves into the $(t-1)$-th row. Thus we need to perform in all $t-2$ straightening
operations, the last one resulting in a multiple of $\la_{t-1}$ as we have to multiply $e_{t-1}$
and $e_{t-1}^{(\la_{t-1} - 1)}$ while moving $e_{t-1}$ into the $(t-1)$-th row. Each straightening
also results in a negative sign. Altogether, the first term after straightening gives
$$(-1)^{t-2} \la_{t-1} \Big(e_1^{(\la_1)} \otimes e_2^{(\la_2)} \otimes  \ldots  e_{t-1}^{(\la_{t - 1})}
\otimes e_{t}^{(\la_t - 1)} \ldots \Big) .$$
By the exact same procedure, the second term, after $t-1$ straightening operations, gives
$$(-1)^{t-1} (\la_{t}-1) \Big(e_1^{(\la_1)} \otimes e_2^{(\la_2)} \otimes  \ldots
e_{t-1}^{(\la_{t - 1})} \otimes e_{t}^{(\la_t - 1)} \ldots \Big) .$$
Altogether we get $(-1)^{t-2} (\la_{t-1} - \la_{t} + 1)$ times a fragment that matches the
corresponding tableau fragment in the canonical tableau $C_\nu$. Evidently the same pattern
will continue as we apply further polarizations. For example after applying the next
polarization, say from the $s$-th row into the $t$-th row, and straightening we will get
an additional multiple of
$(-1)^{s-t-1} (\la_{s-1} - \la_{s} + 1)$
and the resulting tableau fragment with $s$ rows will match the corresponding fragment of
the canonical tableau $C_\nu$. The claimed expression follows after checking easily that
the resulting product of signs matches the claimed sign.

\vskip 0.3cm

Using the work done so far, we can lay out the whole calculation as follows.
$$ \eqalign{ f''(\iota(d'_\nu(C_\nu))) & = f''\left(d'_{\la/1}
\Big(e_1^{(\la_1 - 1)} \otimes e_1 e_2^{(\la_2 - 1)} \otimes \ldots \Big)\right) \cr
\noalign{\vskip 5pt}
& = \sum_{1i_2 \ldots i_k} b_{1 i_2 \ldots i_k} \, g_{i_2 \ldots i_k}''
\Big(e_1^{(\la_1 - 1)} \otimes e_1 e_2^{(\la_2 - 1)} \otimes \ldots \Big) \cr
\noalign{\vskip 5pt}
& = (-1)^{p_1 + \ldots + p_k}  D_\nu \, d'_\nu(C_\nu) ,\cr
}$$
where $D_\nu$ is the appropriate integer obtained using $(*)$ above and the values of
$b_{\underline i}$ obtained near the end of the proof of Lemma C. To finish the proof
we will write an explicit expression for $D_\nu$ and simplify it. First observe that
the product $\prod_{t \in S} (\la_{t-1} - \la_{t} + 1)$ in $(*)$ immediately boils downs
to $a_k \prod_{i_j=1}(a_{j-1}-a_{j}+1)$, where $2 \leq j \leq k$. (The factor $a_k$ is
present because the last row of $\la$, consisting of exactly one box, will always
satisfy the condition defining the set $S$ for all $g_{i_2 \ldots i_k}''$.) So we have
$$ D_\nu = \sum_{\underline i = 1 i_2 \ldots i_k} \left(\prod_{i_j = 0} (h_j - p_j) \right)
\left(a_k \prod_{i_j = 1} (a_{j-1} - a_j +1) \right),$$
where $2 \leq j \leq k$. We will show by induction on the number of blocks $k$ in the
partition $\nu$ that $D_\nu = \ell_1 \ldots \ell_k.$ The case $k=1$ is immediate, since
then the whole expression collapses to just $a_k = a_1 = \ell_1$. When $k > 1$, let $\nu' =$
the partition obtained from $\nu$ by deleting the first block $a_1^{p_1}$. Then the terms
in $D_\nu$ with $i_2=0$ add up to $(h_2 - p_2) D_{\nu'}$, the terms with $i_2 = 1$
add up to  $(a_1 - a_2 +1) D_{\nu'}$ and the terms with any other fixed value of
$i_2$ (i.e., $2 \leq i_2 \leq p_2$) add up to $D_{\nu'}$. Since
$$(h_2 - p_2) + (a_1 - a_2 + 1) + (p_2 - 1) = (a_2 + p_3 + \ldots + p_k) + a_1 - a_2 + p_2
= a_1 + p_2 + \ldots + p_k = \ell_1$$
the proof of Lemma A is complete. This also completes the proof of Theorem 2.1.

\vskip 1cm

\centerline {ACKNOWLEDGMENTS}

\vskip 0.5cm

Theorem 2.1 in its combinatorial form appears in the author's Ph.D. thesis written
under the supervision of David Buchsbaum at Brandeis University. I am grateful
to him for constant encouragement. I also thank Alexei Rudakov and Kari Vilonen
for many discussions. The first version of this paper was written at University
of Massachusetts at Amherst. Finally, I am grateful to the referee for an
exceptionally thorough and detailed review, resulting in numerous improvements
and removal of inaccuracies, especially in the introduction and in Section 1. The
referee's insistence on complete proofs was also instrumental in transforming the
earlier cavalier treatment of lemmas A, B and C into the current form, in particular
the shortened proofs of lemmas A and B.

\vskip 1cm

\centerline  {REFERENCES}

\vskip .5cm

\parindent = 0pt

[Akin] K. Akin, Extensions of symmetric tensors by alternating tensors,
{\it J. Algebra} {\bf 121} (1989), 358--363.

\vskip .2cm

[AB1] K. Akin and D. A. Buchsbaum, Characteristic-free representation
theory of the general linear group,
{\it Adv. in Math.} {\bf 58} (1985), 149--200.

\vskip .2cm

[AB2] K. Akin and D. A. Buchsbaum, Characteristic-free representation
theory of the general linear group II: Homological considerations,
{\it Adv. in Math.} {\bf  72} (1988), 171--210.

\vskip .2cm

[ABW] K. Akin, D. A. Buchsbaum, and J. Weyman, Schur functors and Schur complexes,
{\it Adv. in Math.} {\bf 44} (1982), 207--278.

\vskip .2cm

[Andersen1] H. H. Andersen, Filtrations of cohomology modules for Chevalley groups,
{\it Ann. scient. \'Ec. Norm. Sup.} (4) {\bf 16} (1983), 495--528.

\vskip .2cm

[Andersen2] H. H. Andersen, Filtrations and tilting modules,
{\it Ann. scient. \'Ec. Norm. Sup.} (4) {\bf 30} (1997), 353--366.

\vskip .2cm

[Andersen3] H. H. Andersen, A sum formula for tilting filtrations,
{\it J. Pure Appl. Algebra} {\bf 152} (2000), 17--40.

\vskip .2cm

[Boe] B. D. Boe, A counterexample to the Gabber-Joseph conjecture,
Kazhdan-Lusztig theory and related topics, Chicago 1989,
Contemp. Math. {\bf 139}, Amer. Math. Soc., Providence RI, 1992.

\vskip .2cm

[Boffi] G. Boffi, Characteristic-free decomposition of skew Schur functors,
{\it J. Algebra} {\bf 125} (1989), 288--297.

\vskip .2cm

[BF] D. A. Buchsbaum and D. Flores de Chela, Intertwining numbers: the three-rowed case,
{\it  J. Algebra} {\bf 183} (1996), 605--635.

\vskip .2cm

[CPS]  E. Cline, B. Parshall, and L. Scott, Finite-dimensional algebras and highest weight categories,
{\it J. Reine Angew. Math.} {\bf 391} (1988), 85--99.

\vskip .2cm

[CPSvdK] E. Cline, B. Parshall, L. Scott, and W. van der Kallen,
Rational and generic cohomology, {\it Invent. Math.} {\bf 39} (1977), 143--163.

\vskip .2cm

[CE] A. Cox and K. Erdmann, $\Ext^2$ between Weyl modules for quantum ${\rm GL}_n$,
{\it Math. Proc. Cambridge Philos. Soc.} {\bf 128} (2000), 441--463.

\vskip .2cm

[Donkin1] S. Donkin, On Schur algebras and related algebras I,
{\it J. Algebra} {\bf 104} (1986), 310--328.

\vskip .2cm

[Donkin2] S. Donkin, Skew modules for reductive groups,
{\it J. Algebra} {\bf 113} (1988), 465--479.

\vskip .2cm

[Donkin3] S. Donkin, On Kulkarni's theorems on degree reduction for polynomial modules,
{\it Math. Proc. Cambridge Philos. Soc.} {\bf 134} (2003), 229--237.

\vskip .2cm

[Doty]  S. Doty, The submodule structure of certain Weyl modules for groups of type $A_n$,
{\it J. Algebra} {\bf 95} (1985), 373--383.

\vskip .2cm

[Erdmann] K. Erdmann, $\Ext^1$ for Weyl modules of ${\rm SL}_2(K)$,
{\it Math. Z.} {\bf  218} (1995), 447--459.

\vskip .2cm

[F]  D. Flores de Chela, On intertwining numbers,
{\it J. Algebra} {\bf 171} (1995), 631--653.

\vskip .2cm

[GJ] O. Gabber and A. Joseph, Towards the Kazhdan-Lusztig conjecture,
{\it Ann. scient. \'Ec. Norm. Sup.} (4) {\bf 14} (1981), 261--302.

\vskip .2cm

[Green] J. A. Green, Polynomial representations of ${\rm GL}_n$,
Lecture Notes in Mathematics {\bf 830}, Springer, Berlin-New York, 1980.

\vskip .2cm

[Irving] R. S. Irving, BGG algebras and the BGG reciprocity principle,
{\it J. Algebra} {\bf 135} (1990), 363--380.

\vskip .2cm

[James]  G. D. James, On the decomposition matrices of the symmetric groups III,
{\it J. Algebra} {\bf 71} (1981), 115--122.

\vskip .2cm

[Jantzen] J. C. Jantzen, Representations of algebraic groups, second edition,
Mathematical Surveys and Monographs {\bf 107}, Amer. Math. Soc., 2003.

\vskip .2cm

[Kouwenhoven] F. M. Kouwenhoven, Schur and Weyl functors,
{\it Adv. in Math.} {\bf 90} (1991), 77--113.

\vskip .2cm

[Krop] L. Krop, On the representations of the full matrix semigroup on homogeneous polynomials,
{\it J. Algebra.} {\bf 99} (1986), 370--421.

\vskip .2cm

[Kulkarni1] U. Kulkarni, Skew Weyl modules for ${\rm GL}_n$ and degree reduction for Schur algebras,
{\it J. Algebra} {\bf 224} (2000), 248--262.

\vskip .2cm

[Kulkarni2] U. Kulkarni, A homological interpretation of Jantzen's sum formula, {\it preprint},
arXiv:math.RT/0505371.

\vskip .2cm

[Maliakas] M. Maliakas, Resolutions, homological dimensions, and extensions of
hook representations, {\it Commun. Algebra} {\bf 19} (1991), 2195--2216.

\vskip .2cm

[O-M] M. T. F. Oliveira-Martins, On homomorphisms between Weyl modules for hook partitions,
{\it Linear and Multilinear Algebra} {\bf 23} (1988), 305--323.

\vskip .2cm

[R-G] G. Rondon-Gonzalez, Thesis, Brandeis University, 1998.

\vskip .2cm

[Wen] K. Wen, The composition of intertwining homomorphisms,
{\it Commun. Algebra} {\bf 17} (1989), 587--630.

\end